\documentclass[10pt]{amsart}
\usepackage{latexsym,enumerate}
\usepackage{amsmath,amsthm,amsopn,amstext,amscd,amsfonts,amssymb}
\usepackage{color}
\usepackage{fancybox}
\usepackage{epsfig}
\usepackage{float}
\usepackage{ulem} 

\usepackage[active]{srcltx}

\usepackage{color}

\newcommand{\R}{\mathbb{R}}
\newcommand{\N}{\mathbb{N}}

\newcommand{\car}{{\raise4pt\hbox{$\chi$}}}

\newcommand{\sg}{{\rm \; sign \;}}

\newcommand{\Div}{\hbox{\rm div\,}}

\newcommand{\qin}{\qquad\mbox{ in }\quad}
\newcommand{\qon}{\qquad\mbox{ on }\quad}

\newcommand{\dis }{{\mathcal D}' }
\newcommand{\z }{{\bf z}}
\newcommand{\DM }{\mathcal{DM}^\infty }
\newenvironment{pf}{\noindent{\sc Proof}.\enspace}{\rule{2mm}{2mm}\medskip}

\newtheorem{Theorem}{Theorem}[section]
\newtheorem{Corollary}[Theorem]{Corollary}
\newtheorem{Definition}[Theorem]{Definition}
\newtheorem{Lemma}[Theorem]{Lemma}
\newtheorem{Proposition}[Theorem]{Proposition}

\newtheorem{remark}[Theorem]{Remark}
\newtheorem{Example}[Theorem]{Example}
\newcommand{\res}{\!\!\mathop{\hbox{
                                \vrule height 7pt width .5pt depth 0pt
                                \vrule height .5pt width 6pt depth 0pt}}
                                \nolimits}
\newcommand{\norma}[2]{\|#1\|_{\lower 4pt \hbox{$\scriptstyle #2$}}}

\font\fpe=cmr6

\usepackage{multirow}
\usepackage{array}

\begin{document}

\title[The 1--Laplacian equation with a total variation term]{Elliptic equations involving the 1--Laplacian and a total variation term with $L^{N,\infty}$--data}

\author[M. Latorre and  S. Segura de Le\'on]
{Marta Latorre and  Sergio Segura de Le\'on}

\address{M. Latorre: Departament d'An\`{a}lisi Matem\`atica,
Universitat de Val\`encia,
Dr. Moliner 50, 46100 Burjassot, Spain.
{\it E-mail address:} {\tt marta.latorre@uv.es }}

\address{S. Segura de Le\'on: Departament d'An\`{a}lisi Matem\`atica,
Universitat de Val\`encia,
Dr. Moliner 50, 46100 Burjassot, Spain.
{\it E-mail address:}  {\tt sergio.segura@uv.es }}

\thanks{}


 \bigskip
\begin{abstract}
In this paper we study, in an open bounded set $\Omega\subset\R^N$
with Lipschitz boundary $\partial\Omega$, the Dirichlet problem
for a nonlinear singular elliptic equation involving the
$1$--Laplacian and a total variation term, that is, the
inhomogeneous case of the equation appearing in the level set
formulation of the inverse mean curvature flow.
Our aim is twofold.
On the one hand, we consider data belonging to the Marcinkiewicz space $L^{N,\infty}(\Omega)$, which leads to unbounded solutions.
So, we have to begin introducing the suitable notion of unbounded solution to this problem.
Moreover, examples of explicit solutions are shown.
On the other hand, this equation allows us to deal with many related problems having a different gradient term (see \eqref{prob-prin} below).
It is known that the total variation term induces a regularizing effect on existence, uniqueness and regularity. We focus on analyzing whether those features remain true when general gradient terms are taken. Roughly speaking, the bigger $g$, the better the properties of the solution.
\end{abstract}

\maketitle


\section{Introduction}
In the present paper we deal with the Dirichlet problem for
equations involving the $1$--Laplacian and a total variation term:
\begin{equation}\label{prob-prin}
  \left\{\begin{array}{ll}
  \displaystyle -\Div\Big(\frac{Du}{|Du|}\Big)+g(u)|Du|=f(x)&\hbox{in
   }\Omega\,,\\[3mm]
   u=0 &\hbox{on }\partial\Omega\,,
  \end{array}\right.
\end{equation}
where $\Omega\subset\R^N$ is a bounded open set with Lipschitz
boundary $\partial\Omega$, $g$ stands for a continuous real
function and $f$ is a nonnegative function belonging to the
Marcinkiewicz space $L^{N,\infty}(\Omega)$.

A related class of elliptic problems involving the $p$--Laplacian
operator (defined in $W^{1,p}(\Omega)$ by
 $\displaystyle \Delta_p u=\Div\big(|\nabla u|^{p-2}\nabla
u\big)$, where $p>1$) with a gradient term has been widely
studied. We recall the seminal paper \cite{LL} for a gradient term
of exponent $p-1$ and the systematic study of equations having a
gradient term with natural growth initiated by Boccardo, Murat and
Puel (see \cite{BMP1, BMP2, BMP3}). The variational approach
searches for solutions in the Sobolev space $W_0^{1,p}(\Omega)$ and
considers data belonging to its dual $W^{-1,p'}(\Omega)$.  (In the
setting of Lebesgue spaces, data are naturally taken in
$L^{\frac{Np}{Np-N+p}}(\Omega)$ as a consequence of the Sobolev
embedding.)

We point out that the natural space to look for a solution to
problem \eqref{prob-prin} should be the Sobolev space
$W_0^{1,1}(\Omega)$ and the space of data, from a variational point
of view, should be its dual $W^{-1,\infty}(\Omega)$. The Sobolev
embedding Theorem and duality arguments lead to consider as {\it
the right} function space of data the space $L^N(\Omega)$ (among
the Lebesgue spaces) and $L^{N,\infty}(\Omega)$ (among the Lorentz
spaces). Evidences that the norm of $L^{N,\infty}(\Omega)$ is
suitable enough to deal with this kind of problems can be found in
\cite{CT, MST1}. As far as the energy space is concerned, we
cannot search for solutions in $W_0^{1,1}(\Omega)$, which is not
reflexive, and we have to extend our setting to the larger space
$BV(\Omega)$, the space of all functions of bounded variation.
Therefore, our framework is the following: given a nonnegative
$f\in L^{N,\infty}(\Omega)$, find $u\in BV(\Omega)$ that solves
problem \eqref{prob-prin} in an appropriate sense which will be introduced
below (see Definition \ref{def1}).

Two important cases of problem \eqref{prob-prin} have already been
studied. When $g(s)\equiv 0$ we obtain just the $1$--Laplacian
operator: $ \displaystyle -\Div\Big(\frac{Du}{|Du|}\Big)$. There is
a big amount of literature on this equation in recent years,
starting in \cite{K}. Other papers dealing with this equation are
\cite{ABCM, BCN, CT, D, KF, MST1}. The interest in studying such a
case came from an optimal design problem in the theory of torsion
and related geometrical problems (see \cite{K}) and from the
variational approach to image restoration (see \cite{ABCM} and
also \cite{ACM} for a review on the development of variational
models in image processing). The suitable concept of solution to
handle the Dirichlet problem for this kind of equations was
introduced in \cite{ABCM}. In this paper, a meaning for the
quotient $\displaystyle \frac{Du}{|Du|}$ (involving Radon
measures) is given through a vector field $\z\in
L^\infty(\Omega;\R^N)$ satisfying $\|\z\|_\infty\le1$ and $(\z,
Du)=|Du|$ as measures. This vector field also gives sense to the
boundary condition in a weak sense. The meaning of all expressions
in which appear vector fields relies on the theory of
$L^\infty$--divergence--measure fields (see \cite{An} and \cite{CF}).

On the other hand, when $g(s)\equiv 1$, we get  $\displaystyle
-\Div\Big(\frac{Du}{|Du|}\Big)+|Du|$, which occurs in the level
set formulation of the inverse mean curvature flow (see
\cite{Huisken}, related developments can be found in
\cite{Huisken2, Moser, Moser2}).
 The framework of these papers,
however, is different since $\Omega$ is unbounded. Furthermore,
the concept of solution is based on the minimization of certain
functional and does not coincide with which has been considered in the
previous case. This operator has also been studied in a bounded
domain in \cite{MS}, where it is proved the existence and
uniqueness of a bounded solution for a datum regular enough.

It is worth noting that, contrary to what happens in the
$p$--Laplacian setting with $p> 1$, features of solutions to
problem \eqref{prob-prin} with $g(s)\equiv 0$ are very different
to those with $g(s)\equiv 1$. Indeed, the presence of the gradient
term has a strong regularizing effect because in the first case the
following facts hold:
\begin{enumerate}
    \item[(i)]  Existence of $BV$--solutions is only guaranteed for data
small enough, for large data solutions become infinity in a set of
positive measure.
    \item[(ii)] There is no uniqueness at all: given a solution $u$,
    we also obtain that $h(u)$ is a solution, for every smooth increasing function $h$.
 \end{enumerate}
Whereas, in the second case, the properties are:
\begin{itemize}
    \item[(i)]  There is always a solution, even in the case where the datum is large.
    \item[(ii)] An uniqueness result holds.
\end{itemize}
Regarding regularity of solutions, even an equation related to the case $g(s)\equiv 0$ like $u-\Div\Big(\frac{Du}{|Du|}\Big)=f(x)$ (for which existence and uniqueness hold) has solutions with jump part. On the contrary, solutions to problem \eqref{prob-prin} with $g(s)\equiv 1$ have no jump part.
Moreover, solutions to $u-\Div\Big(\frac{Du}{|Du|}\Big)=f(x)$ satisfy
the boundary condition only in a weak sense (and in
general, $u|_{\partial\Omega}\ne 0$), while if $g(s)\equiv 1$,
then  the boundary condition holds in the trace sense, that is,
the value is attained pointwise on the boundary.

We point out that the situation concerning existence is rather
similar to that shown in studying problem
\begin{equation}\label{hardy}
  \left\{\begin{array}{ll}
  \displaystyle -\Delta u+|\nabla u|^2=\lambda \frac{u}{|x|^2}&\hbox{in
   }\Omega\,,\\[3mm]
   u=0 &\hbox{on }\partial\Omega\,,
  \end{array}\right.
\end{equation}
in domains satisfying $0\in\Omega$, since the presence of the
quadratic gradient term induces a regularizing effect (see
\cite{APP} and \cite{ABPP}, see also Remark \ref{nota} below). Indeed, existence of a positive
solution to \eqref{hardy} can be proved for all $\lambda>0$, while
if the gradient term does not appear, solutions can be
expected only for $\lambda$ small enough, due to Hardy's
inequality.

Our purpose is to study the role of the function $g$ on the above
features satisfied by the solutions. Roughly speaking, we see that
the bigger $g$, the better the properties of the solution. The
standard case occurs when $g(s)\ge m>0$ for all $s\ge0$ and the
situation degenerates as soon as $g(s)$ touch the $s$--axis.

We begin by considering the case $g(s)=1$ for all $s\ge0$. To get an
idea of the difficulties one finds, let us recall previous works
on this subject. As mentioned, this problem was already handled in
\cite{MS} for data $f\in L^q(\Omega)$, with $q>N$. This condition
is somewhat artificial and was taken in this way due to the
necessity of obtaining bounded solutions. This necessity derives
from the use of the theory of $L^\infty$--divergence--measure
fields. It was initiated in \cite{An}, where a sense is provided
with the dot product $(\z, Du)$, where $\z\in L^\infty(\Omega; \R^N)$
satisfies that $\Div\z$ is a Radon measure and $u\in
BV(\Omega)\cap L^\infty(\Omega)$ is a continuous function. In a
different way, it was later developed in \cite{CF} for a -possibly
discontinuous- function $u\in BV(\Omega)\cap L^\infty(\Omega)$ (see also
\cite{C, MST2} for a point of view closest to that of \cite{An}).
Since we must expect unbounded solutions starting from the most
natural space of data $L^{N,\infty}(\Omega)$, the first result we
need is to give sense to the dot product $(\z, Du)$ when $u \in
BV(\Omega)$ can be unbounded. This was achieved in \cite{ADS}, but we
include it for the sake of completeness.

 Endowed with this tool, in the first part of this paper, we prove an existence and
uniqueness result for problem \eqref{prob-prin} in the particular
case $g(s)\equiv 1$. The second part is fully devoted to our main
concern, that is, to search for the properties that solutions to
problem \eqref{prob-prin} satisfy for different functions $g$. For
better understanding, we summarize the results we will see in the
table below.

\vskip7mm

\centerline{\fpe
\begin{tabular}{|c||c|c|c|}
 \hline
 &&&\\[-0.05cm]
 Function $g(s)$ & Existence & Uniqueness & Regularity \\[0.4cm]
 \hline
 \hline
 &&&\\[-0.1cm]
 \multirow{2}*{$0<m\le g(s)$} & \multirow{2}*{For every datum$^{(1)}$} & \multirow{2}*{Yes$^{(1)}$} & No jump part$^{(1)}$\\
 &&& Better summability$^{(2)}$ \\[0.2cm]
 \hline
 &&&\\[-0.1cm]
$g$ vanishes at some points & \multirow{2}*{For every datum$^{(3)}$} & \multirow{2}*{Yes$^{(3)}$} & \multirow{2}*{No jump part$^{(3)}$} \\
$g\notin L^1([0,\infty[)$ &  &  & \\[0.25cm]
\hline
&&&\\[-0.1cm]
$g$ vanishes at infinity & For every datum$^{(4)}$, with & \multirow{2}*{Yes$^{(4)}$} & \multirow{2}*{No jump part$^{(4)}$}\\
$g\notin L^1([0,\infty[)$ &  another concept of solution$^{(5)}$ &  & \\[0.2cm]
\hline
&&&\\[-0.1cm]
 $g\in L^1([0,\infty[)$ & For data small enough$^{(6,7)}$ & Yes$^{(7)}$ & No jump part$^{(7)}$ \\[0.25cm]
 \hline
 &&&\\[-0.1cm]
\multirow{2}*{$g$ vanishes on an interval} & \multirow{2}*{For data small enough$^{(8)}$} & \multirow{2}*{No$^{(9)}$} & With jump part$^{(10)}$\\
&&&  No boundary condition$^{(11)}$\\[0.2cm]
\hline
\end{tabular}}

\smallskip

{\fpe{ (1) Theorem \ref{teounicidad} and Theorem \ref{teoalpha} \qquad (2) Proposition \ref{propq} \qquad (3) Theorem \ref{teo-finitos-ceros}\qquad (4) Theorem \ref{teoproblemag}

(5) Definition \ref{nueva} and Example \ref{toca-infin}\qquad(6) Example \ref{no-sol-radial}
 \qquad (7) Theorem \ref{raro1} \qquad (8) Remark \ref{no-unicidad1}

(9) Remark \ref{no-unicidad1} and Remark \ref{no-unicidad2} \qquad (10) Example \ref{ejemplo-disc} \qquad (11) Example \ref{ejemplo-no-frontera} }}

\vskip7mm

The plan of this paper is the following. Section 2 is dedicated to
preliminaries, we introduce our notation and some properties of
the spaces $BV(\Omega)$ and $L^{N, \infty}(\Omega)$. In Section 3
we generalize the theory of $L^\infty$--divergence--measure fields
to take pairings $(\z, Du)$ of a certain vector field $\z$ and any
$u\in BV(\Omega)$. This theory is applied in Section 4 to extend
the result of existence and uniqueness of \cite{MS} to $L^{N,
\infty}(\Omega)$--data. In Section 5 we show explicit radial
examples of solutions. Section 6 is devoted to study the standard
cases of problem \eqref{prob-prin}, those where $g(s)$ is bounded
from below by a positive constant. A non standard case is shown in
Section 7 with $g(s)$ touching the $s$--axis; in this case we need
to change our definition of solution since solutions no longer
belong to $BV(\Omega)$. Finally, in Section 8 we deal with really
odd cases for which the considered properties are not necessarily
satisfied.

\section*{Acknowledgements}
This research has been partially supported  by the Spanish Mi\-nis\-te\-rio de Econom\'{\i}a y Competitividad and
FEDER, under project MTM2015--70227--P.
The first author was also supported by Ministerio de Econom\'{\i}a y
Competitividad under grant BES--2013--066655.

The authors would like to thank Salvador Moll for bringing paper \cite{GMP} to our attention.

\section{Preliminaries}

In this Section we will introduce some notation and auxiliary results which will be used
throughout this paper. In what follows, we will consider $N\ge2$,
and $\mathcal H^{N-1}(E)$ will denote the $(N - 1)$--dimensional
Hausdorff measure of a set $E$ and $|E|$ its
Lebesgue measure.

  In this paper, $\Omega$ will always denote an open subset of
  $\R^N$ with Lipschitz boundary. Thus, an outward normal unit
  vector $\nu(x)$ is defined for $\mathcal H^{N-1}$--almost every
  $x\in\partial\Omega$.
 We will make use of the usual Lebesgue and Sobolev
 spaces, denoted by $L^q(\Omega)$  and $W_0^{1,p}(\Omega)$,
 respectively.

We recall that for a Radon measure $\mu$ in $\Omega$ and a Borel set
$A\subseteq\Omega$ the measure $\mu\res A$ is defined by $(\mu\res A)(B)=\mu(A\cap B)$
for any Borel set $B\subseteq\Omega$. If a measure $\mu$ is such that $\mu
 = \mu \res A$ for a certain Borel set $A$, the measure $\mu$ is
 said to be concentrated on $A$.

The truncation function will be use throughout this paper. Given
$k>0$, it is defined by
\begin{equation}\label{trun}
    T_k(s)=\min\{|s|, k\}\sg (s)\,,
\end{equation} for all $s\in\R$. Moreover, we define another auxiliary real function by
\begin{equation}\label{G-k-fun}
    G_k(s)=\big(s-T_k (s)\big)\sg (s)\,.
\end{equation}

\subsection{The energy space}

The space of all functions of finite variation, that is the space
of those $u\in L^1(\Omega)$ whose distributional gradient is a Radon
measure with finite total variation, is denoted by $BV(\Omega)$.
This is the natural energy space to study the problems we are interested
in. It
is endowed with the norm defined by
 $$ \|u\|=\int_\Omega |u|\, dx+ \int_\Omega|Du|\,,$$
 for any $u\in BV(\Omega)$. An equivalent norm, which we will use in the sequel, is given by
 $$\displaystyle \|u\|_{BV(\Omega)}=\int_{\partial\Omega}
|u|\, d\mathcal H^{N-1}+ \int_\Omega|Du|\,.$$

For every
$u \in BV(\Omega)$, the Radon measure $Du$ is decomposed into its
absolutely continuous and singular parts with respect to the
Lebesgue measure: $Du = D^a u + D^s u$.
We denote by $S_u$ the set of all $x\in\Omega$ such that $x$ is
not a Lebesgue point of $u$, that is, $x\in\Omega\backslash S_u$ if there exists $\tilde{u}(x)$ such that
$$\lim_{\rho \downarrow 0} \frac{1}{|B_{\rho}(x)|}
\int_{B_{\rho}(x)} \vert u(y) - \tilde{u}(x) \vert \, dy = 0\,.$$
 We say that $x \in \Omega$ is an {\it
approximate jump point of } $u$ if there exist two real numbers $u^+(x)
>u^-(x) $ and $\nu_u(x) \in S^{N-1}$ such that
$$\lim_{\rho \downarrow 0} \frac{1}{|B_{\rho}^+(x,\nu_u(x))|}
\int_{B_{\rho}^+(x,\nu_u(x))} \vert u(y) - u^+(x) \vert \, dy = 0\,,$$
$$\lim_{\rho \downarrow 0} \frac{1}{|B_{\rho}^-(x,\nu_u(x))|}
\int_{B_{\rho}^-(x,\nu_u(x))} \vert u(y) - u^-(x) \vert \, dy = 0\,,$$ where
$$B_{\rho}^+(x,\nu_u(x)) = \{ y \in B_{\rho}(x) \ : \ \langle y - x, \nu_u(x) \rangle >
0 \} $$ and
$$B_{\rho}^-(x,\nu_u(x)) = \{ y \in B_{\rho}(x) \ : \ \langle y - x, \nu_u(x) \rangle <
0 \}\,.$$
 We denote by $J_u$ the set of all approximate jump points of
$u$. By the Federer--Vol'pert Theorem \cite[Theorem 3.78]{AFP}, we know that $S_u$ is countably $\mathcal H^{N-1}$--rectifiable and $\mathcal H^{N-1}(S_u \backslash J_u) = 0$. Moreover, $Du \res J_u = (u^+ - u^-) \nu_u \mathcal H^{N-1} \res J_u$. Using
$S_u$ and $J_u$, we may split $D^su$ in two parts: the {\it jump}
part $D^j u$ and the {\it Cantor} part $D^c u$ defined by
$$D^ju = D^su \res J_u  \ \ \ {\rm and} \ \ D^c u = D^su \res (\Omega \backslash S_u)\,.$$
Then, we have
$$D^j u = (u^+ - u^-) \nu_u \mathcal H^{N-1} \res J_u\,.$$
Moreover, if $x \in J_u$, then $\nu_u(x) = \frac{Du}{| D u |}(x)$ and
$\frac{Du}{| D u |}$ is the Radon--Nikod\'ym derivative of $Du$
with respect to its total variation $| D u |$.

The precise representative
$u^* : \Omega \backslash(S_u \backslash J_u) \rightarrow \R$ of
$u$ is defined as equal to $\tilde{u}$ on $\Omega \backslash S_u$
and equal to $\frac{u^- + u^+}{2}$ on $J_u$. It is well known (see
for instance \cite[Corollary 3.80]{AFP}) that if $\rho$ is a
symmetric mollifier, then the mollified functions $u \star
\rho_{\epsilon}$ pointwise converges to $u^*$ in its domain.

  A compactness result in $BV(\Omega)$ will be used several times in what follows.
  It states that every sequence that is bounded in $BV(\Omega)$ has a
subsequence which strongly converges in $L^1(\Omega)$ to a certain
$u\in BV(\Omega)$ and the subsequence of gradients $*$--weakly
converges to $Du$ in the sense of measures.

 To pass to the limit we will often apply that some functionals defined on $BV(\Omega)$ are
 lower semicontinuous with respect to the convergence in $L^1(\Omega)$.
  The most important are
  the functionals
  defined by
  \begin{equation}\label{semc}
 u\mapsto\int_{\Omega}|Du|
 \end{equation}
 and
   \begin{equation}\label{semcon}
  u\mapsto\int_\Omega|Du|+\int_{\partial\Omega}|u|\,d\mathcal
  H^{N-1}\,.
  \end{equation}
  In the same way,
  it yields that each $\varphi\in C_0^1(\Omega)$ with
  $\varphi\ge0$ defines a functional
  $$
  u\mapsto\int_{\Omega}\varphi\,|Du|\,,
  $$
  which is lower semicontinuous in $L^1(\Omega)$.

  Finally, we recall that the notion of trace can be extended to any $u\in BV(\Omega)$ and this fact allows us to interpret it as the boundary values of $u$ and to write $u\big|_{\partial \Omega}$. Moreover, it holds that the trace is a linear bounded operator $BV(\Omega)\to  L^1(\partial\Omega)$ which is onto.

For further information on functions of bounded variation, we refer to \cite{AFP, EG, Zi}.

\subsection{The data space}
Given a measurable function $u:\Omega\to\R$, we denote by $\mu_u$
the distribution function of $u$: the function $\mu_u: [0,
+\infty[\rightarrow [0, +\infty[$ defined by
$$\mu_u(t)=|\{x\in\Omega :|u(x)|>t\}|\,,\quad t\ge0\,.$$

 For $1< q<\infty$, the  space $L^{q,\infty}(\Omega)$, known as
Marcinkiewicz or weak-Lebesgue space, is the space of Lebesgue
measurable functions $u:\Omega\to\R$ such that
 \begin{equation}\label{quasi-norm}
[u]_q=\sup_{~t>0}t\,\mu_u(t)^{1/q}<+\infty\,.
\end{equation}
The relationship with Lebesgue spaces is
given by the following inclusions
 $$L^q(\Omega) \hookrightarrow  L^{q,\infty}(\Omega) \hookrightarrow
L^{q-\varepsilon}(\Omega)\,,$$
for suitable $\varepsilon >0$.
We point out that expression \eqref{quasi-norm} defines a quasi--norm which is not a norm in
$L^{q,\infty}(\Omega)$. (For a suitable norm in this space see \eqref{m1-norm}, \eqref{m2-norm} and \eqref{m3-norm} below).

  Some properties of Lorentz spaces $L^{q,1}(\Omega)$
(with $1<q<\infty$) must be applied throughout this paper. To
begin with, we define the decreasing rearrangement of $u$ as the
function $u^\star:\,]0,|\Omega|]\to\R^+$ given by
$$u^\star(s)=\sup\{t>0:\mu_u(t)>s\}\,, \qquad s\in \, ]0,|\Omega|]\,,$$
(the main properties of rearrangements can be found in \cite{BS, H, Zi}).
In terms of $u^\star$, the quasi-norm \eqref{quasi-norm} becomes
\begin{equation}\label{m-norm}
  [u]_q=\sup_{s>0}\,\{s^{1/q}u^\star(s)\}\,.
\end{equation}

We say that a measurable function $u\>:\>\Omega\to\R$
belongs to $L^{q,1}(\Omega)$ if
\begin{equation}\label{l-norm}
  \|u\|_{L^{q,1}(\Omega)}=\frac1q\int_0^\infty s^{1/q}u^\star(s)\frac{ds}s
\end{equation}
is finite. This expression defines a norm (see \cite[Theorem
5.13]{BS}). The classical paper where these spaces are
systematically studied is \cite{H} (see also \cite{BS, Zi}). Some
important properties of Lorentz spaces are:
\begin{enumerate}
\item $L^{q,1}(\Omega)$ is a Banach space endowed with the norm defined by \eqref{l-norm}.
\item Simple functions are dense in $L^{q,1}(\Omega)$.
\item The norm \eqref{l-norm} is absolutely continuous.
\end{enumerate}

 Concerning duality, the
Marcinkiewicz space $L^{q^\prime,\infty}(\Omega)$ is the dual space of $L^{q,1}(\Omega)$.
 Indeed, it follows from a Hardy--Littlewood inequality that if $f\in L^{q^\prime,\infty}(\Omega)$ and $u\in L^{q,1}(\Omega)$, then $fu\in L^1(\Omega)$ and a H\"older type inequality holds:
\begin{multline*}
  \Big|\int_\Omega fu\, dx\Big|\le \int_0^\infty f^\star(s) u^\star(s)\, ds=\int_0^\infty s^{1/q^\prime}f^\star(s) s^{1/q}u^\star(s)\, \frac{ds}s \\
  \le
  q [f]_{q^\prime}\|u\|_{L^{q,1}(\Omega)}\,.
\end{multline*}
Thus,
\begin{equation}\label{m1-norm}
  \|f\|_{L^{q^\prime,\infty}(\Omega)}=\sup \left \{ \frac{\Big|\int_\Omega fu\, dx\Big|}{\|u\|_{L^{q,1}(\Omega)}}\> : \>
  u\in L^{q,1}(\Omega)\backslash \{0\} \right \}
\end{equation}
defines a norm in the Marcinkiewicz space and $\|f\|_{L^{q^\prime,\infty}(\Omega)}\le q\, [f]_{q^\prime}$ holds.
Taking into account that if $E\subset \Omega$ is a measurable set of positive measure and $u=|E|^{-\frac1q}\chi_E$,
then $\|u\|_{L^{q,1}(\Omega)}=1$ and also applying the density of simple functions, we deduce that
\begin{multline}\label{m2-norm}
\|f\|_{L^{q^\prime,\infty}(\Omega)}=\sup\left\{\Big|\int_\Omega fu\, dx\Big|\>:\> u=|E|^{-\frac1q}\chi_E\,, \hbox{ with } |E|>0\right\}\\
=\sup\left\{|E|^{-1/q}\int_E |f|\, dx\>:\> |E|>0\right\}\,.
\end{multline}
This implies $[f]_{q^\prime}\le \|f\|_{L^{q^\prime,\infty}(\Omega)}$, so that, the quasi--norm $[\, \cdot\, ]_{q^\prime}$ is equivalent to the norm
$\|\cdot\|_{L^{q^\prime,\infty}(\Omega)}$.
It also yields
\begin{equation}\label{m3-norm}
  \|f\|_{L^{q^\prime,\infty}(\Omega)}=\sup_{s>0}\,\{s^{1/q^\prime}f^{\star\star}(s)\}\,,
\end{equation}
where $\displaystyle f^{\star\star}(s)=\frac1s\int_0^sf^\star(\sigma)\, d\sigma\,.$

On the other hand, we recall that Sobolev's inequality can be
improved in the context of Lorentz spaces (see \cite{Al}): the
continuous embedding
 \begin{equation}\label{inclusion1}
W_0^{1,1}(\Omega)\hookrightarrow L^{\frac N{N-1},1}(\Omega)
 \end{equation}
 holds.
  The best constant in this embedding will be denoted as
 \begin{equation}\label{s-cons}
  S_N=\sup\bigg\{\frac{\|u\|_{L^{\frac{N}{N-1},1}(\Omega)}}{\int_\Omega|\nabla u|\, dx}\>:\>
 u\in W_0^{1,1}(\Omega)\backslash\{0\}\bigg\}\,.
  \end{equation}
Its value is  known:
   \begin{equation}\label{s-value}
   S_N=\frac{\Gamma\big(\frac N2+1\big)^{1/N}}{N\sqrt \pi}=\frac1{NC_N^{1/N}}\,,
\end{equation}
where $C_N$ denotes the
measure of the unit ball in $\R^N$.
(We explicitly point out that this is the value for the best constant having in mind the norm in the Lorentz space as defined in \eqref{l-norm}.)
Furthermore, by an approximation argument, this inclusion may be
extended to BV--functions with the same best constant $S_N$ (see, for instance, \cite{Zi}):
 \begin{equation}\label{inclusion2}
BV(\Omega)\hookrightarrow L^{\frac N{N-1},1}(\Omega)\,.
 \end{equation}
 It is worth remarking that the supremum in \eqref{s-cons} is attained in $BV(\Omega)$.

 As a consequence of this embedding, given $f\in L^{N,\infty}(\Omega)$
 and $u\in BV(\Omega)$, it yields $fu\in L^1(\Omega)$.
This fact will be essential in what follows.

  Another fact concerning Lorentz spaces and duality is in order.
  We will denote by $W^{-1, q^\prime}(\Omega)$ the dual space of
$W^{1,q}_0(\Omega)$, $1\le q<\infty$.
 Here we recall just that the norm in $W^{-1,\infty}(\Omega)$
is given by
 \begin{equation}\label{novnorm}
 \|\mu\|_{W^{-1,\infty}(\Omega)}=\sup \left \{ \big|< \mu, u>_{W^{-1,\infty}(\Omega),W_0^{1,1}(\Omega)}\big|\> : \>
\int_{\Omega} |\nabla u| \, dx \le 1 \right \}\,.
 \end{equation}
  Since the norm in $L^{\frac N{N-1},1}(\Omega)$ is absolutely continuous, it follows that $C_0^\infty(\Omega)$ is dense in $L^{\frac N{N-1},1}(\Omega)$.
  A duality argument shows that $L^{N,\infty}(\Omega)\hookrightarrow W^{-1,\infty}(\Omega)$ and, having in mind \eqref{m1-norm} and \eqref{s-cons}, we obtain: if $f\in L^{N,\infty}(\Omega)$, then
\begin{multline*}
\|f\|_{L^{N,\infty}(\Omega)}=\sup \left \{ \frac{\Big|\int_\Omega
fu\, dx\Big|}{\|u\|_{L^{\frac N{N-1},1}(\Omega)}}\> : \>
  u\in W_0^{1,1}(\Omega)\backslash \{0\} \right \}\\
=\sup \left \{ \frac{\Big|\int_\Omega fu\, dx\Big|}{\int_\Omega |\nabla
u|\,dx}\cdot \frac{\int_\Omega |\nabla u|\, dx}{\|u\|_{L^{\frac
N{N-1},1}(\Omega)}}\> : \>
  u\in W_0^{1,1}(\Omega)\backslash \{0\} \right \}\ge S_N^{-1} \|f\|_{W^{-1,\infty}(\Omega)}\,.
\end{multline*}
Therefore,
\begin{equation}\label{rel-norm}
\|f\|_{W^{-1,\infty}(\Omega)}\le \frac1{N
C_N^{1/N}}\|f\|_{L^{N,\infty}(\Omega)}\,,
\end{equation}
for every $f\in L^{N,\infty}(\Omega)$.
(For a related equality in a ball, see \cite[Remark 3.3]{MST1}).

\section{Extending Anzellotti's theory}

In this section we will study some properties involving
divergence--measure vector fields and functions of bounded
variation. Our aim is to extend the Anzellotti theory.

 Following \cite{CF} we define $\DM(\Omega)$ as the space of all
 vector fields $\z\in L^\infty(\Omega;\R^N)$ whose divergence
 in the sense of distributions is a Radon measure with finite total variation, i.e., $\z\in \DM(\Omega)$ if and only if
 $\hbox{div\,} \z$ is a Radon measure belonging to $W^{-1,\infty}(\Omega)$.

The theory of $L^\infty$--divergence--measure vector fields is due
to G. Anzellotti \cite{An} and, independently, to G.--Q. Chen and
H. Frid \cite{CF}. In spite of their different points of view,
both approaches introduce the normal trace of a vector field
through the boundary and establish the same generalized
Gauss--Green formula. Both two also define the pairing $(\z,Du)$ as a
Radon measure
where $\z\in\DM(\Omega)$ and $u$ is a certain $BV$--function. However, they differ in handling this concept.
While in \cite{An} it is only considered continuous functions belonging to
 $BV(\Omega)\cap L^\infty(\Omega)$ and the inequality
  \begin{equation}\label{cruc}
  |(\z,Du)|\le \|\z\|_\infty|Du|
  \end{equation}
  is proved for those functions; in \cite{CF}, general $u\in BV(\Omega)\cap L^\infty(\Omega)$ are
  considered but it is only shown that the Radon measure $(\z,Du)$ is
  absolutely continuous with respect to $|Du|$. In the present
  paper we need that the inequality \eqref{cruc} holds for
  every $u\in BV(\Omega)$ and every $\z\in\DM(\Omega)$ satisfying a certain condition (see
  Corollary \ref{clave} below). That is why the way by which the pairings $(\z,Du)$
  are obtained will be essential in our work. This is the reason for extending
   the Anzellotti approach in this Section.

  We finally point out that the theory
  of divergence--measure fields has been extended later (see \cite{CTZ} and \cite{Zi2}).

  We begin by recalling a result proved in \cite{CF}.

   \begin{Proposition}\label{absolcont}
  For every $\z\in\DM(\Omega)$,
  the measure $\mu=\hbox{\rm div\,}\z$ is absolutely continuous with respect to
 $\mathcal H^{N-1}$, that is,  $|\mu|\ll\mathcal H^{N-1}$.
 \end{Proposition}

 Consider now $\mu=\hbox{\rm div\,}\z$ with $\z\in\DM(\Omega)$ and let $u\in
BV(\Omega)$; then the precise representative $u^*$ of  $u$ is
equal $\mathcal H^{N-1}$--a.e. to a Borel function; that is, to
$\lim_{\varepsilon\to0}\rho_\varepsilon\star u$, where $(\rho_\varepsilon)$
is a symmetric mollifier. Then, it is deduced from the previous
Proposition that $u^*$ is equal $\mu$--a.e. to a Borel function.
So, given $u\in BV(\Omega)$, its precise representative $u^*$ is always $\mu$--measurable.
Moreover, $u\in BV(\Omega)\cap L^\infty(\Omega)$ implies
 $u\in L^\infty (\Omega,\mu)\subset L^1 (\Omega,\mu)$.

\subsection{Preservation of the norm}

  We point out that every $\Div\z$,
  with $\z\in \DM(\Omega)$, defines a functional on
  $W_0^{1,1}(\Omega)$ by
 \begin{equation}\label{def:1}
  \langle\Div\z, u
  \rangle_{W^{-1,\infty}(\Omega),W_0^{1,1}(\Omega)}=-\int_\Omega\z\cdot\nabla
  u\,dx\,.
  \end{equation}
  To express this functional in terms of an integral with respect to the measure $\mu=\hbox{\rm div\,}\z$,
  we need the following Meyers--Serrin type theorem (see \cite[Theorem 3.9]{AFP} for its extension to $BV$--functions).

     \begin{Proposition}\label{M-S}
Let $\mu = \Div \z$, with $\z\in\DM(\Omega)$.
For every $u \in BV (\Omega)\cap L^\infty(\Omega)$ there exists a sequence
$(u_n)_n$ in $W^{1,1}(\Omega) \cap C^\infty(\Omega) \cap L^\infty(\Omega)$ such that
  \begin{equation*}
  \begin{array}{ll}
  \hbox{\rm (1) } u_n\to u^*\quad\hbox{ in } L^1(\Omega,\mu)\,.\\
  \\
  \hbox{\rm (2) } \int_\Omega|\nabla u_n|\,dx\to|Du|(\Omega)\,.\\
  \\
  \hbox{\rm (3) } u_n|_{\partial\Omega}=u|_{\partial\Omega}\hbox{ for all }
  n\in\N\,.\\
  \\
  \hbox{\rm (4) }
  |u_n(x)|\le\|u\|_\infty\ \ |\mu|\hbox{--a.e. for all }
  n\in\N\,.
  \end{array}
  \end{equation*}
Moreover, if $u \in W^{1,1}(\Omega) \cap L^\infty(\Omega)$,
then one may find $u_n$ satisfying, instead of (2), the condition
$$\hbox{\rm (2') } u_n\to u \hbox{ in } W^{1,1}(\Omega)\,.$$
 \end{Proposition}

Since  $$
  -\int_\Omega\z\cdot\nabla
  \varphi\,dx=\int_\Omega \varphi\, d\mu
  $$
  holds for every $\varphi\in C_0^\infty(\Omega)$, it is easy to
  obtain this equality for every $W_0^{1,1}(\Omega) \cap
  C^\infty(\Omega)$.
  Given $u\in W_0^{1,1}(\Omega)\cap L^\infty(\Omega)$ and applying Proposition \ref{M-S}, we may find a sequence $(u_n)_n$ in $W_0^{1,1}(\Omega) \cap C^\infty(\Omega)$ satisfying (1) and (2').
   Letting $n$ go to infinity, it follows from
  $$
  -\int_\Omega\z\cdot\nabla
  u_n\,dx=\int_\Omega u_n\, d\mu
  $$
  for every $n\in \N$, that
  $$
  -\int_\Omega\z\cdot\nabla
  u\,dx=\int_\Omega u^*\, d\mu
  $$
  and so
  $$
  \langle\Div\z, u
  \rangle_{W^{-1,\infty}(\Omega),W_0^{1,1}(\Omega)}=\int_\Omega u^*\,
  d\mu
  $$
  holds for every $u\in W_0^{1,1}(\Omega)\cap L^\infty(\Omega)$. Then the
  norm of this functional is given by
  $$
  \|\mu\|_{W^{-1,\infty}(\Omega)}=\sup\left\{\Big|\int_\Omega u^*\, d\mu\Big|\>:\, u\in W^{1,1}_0(\Omega)\cap L^{\infty}(\Omega) ,\>\text{ with }\>\|u\|_{W^{1,1}_0(\Omega)}\le
  1\right\}\,.
  $$
where $\|u\|_{W^{1,1}_0}=\int_\Omega|\nabla u|\, dx$.
We have seen that $\mu=\hbox{\rm div\,}\z$ can be extended from
  $W_0^{1,1}(\Omega)$ to $BV(\Omega)\cap L^\infty(\Omega)$.
Next, we will prove that this extension can be given as an integral with respect to $\mu$ and it preserves the
  norm. To this end, the following Lemma, stated in \cite{An},
  will be applied.

  \begin{Lemma}
 For every $u\in BV(\Omega)$
 --so that $u\big|_{\partial\Omega}\in L^1(\partial\Omega)$--,
 there exists a sequence $(w_n)_n$ in
 $W^{1, 1}(\Omega)\cap C(\Omega)$ such that
 \begin{equation*}\begin{array}{ll}
 \hbox{\rm (1) }
 w_n|_{\partial\Omega}=u|_{\partial\Omega}\,.\\
\\
 \hbox{\rm (2) } \displaystyle\int_\Omega|\nabla
 w_n|\,dx\le\displaystyle\int_{\partial\Omega}|u|\,d\mathcal H^{N-1}+\frac1n\,.\\
\\
 \hbox{\rm (3) } \displaystyle\int_\Omega| w_n|\,dx\le\frac1n\,.\\
\\
 \hbox{\rm (4) } w_n(x)=0  \quad\hbox{ if }\quad\hbox{\rm dist}(x,
 \partial\Omega)>\frac1n\,.\\
\\
 \hbox{\rm (5) } w_n(x)\to0  \quad\hbox{ for all }\quad x
 \in\Omega\,.\hskip5cm
 \end{array}\end{equation*}

 Moreover, if $u\in BV(\Omega)\cap L^\infty(\Omega)$, then
 $w_n\in L^\infty(\Omega)$ and $\|w_n\|_\infty\le\|u{\big|_{\partial\Omega}}\|_\infty$ for
 all $n\in\N$.
  \end{Lemma}

 \begin{Theorem}\label{ext}
 Let $\z\in \DM(\Omega)$ and denote $\mu=\Div\z$. Then, the functional given by \eqref{def:1} can be extended to $BV(\Omega) \cap L^{\infty}(\Omega)$ as an integral with respect to $\mu$ and its norm satisfies
 \[
 \|\mu\|_{W^{-1,\infty}(\Omega)} = \sup\left\{\Big|\int_{\Omega} u^*\,d\mu\Big| \; : \; u\in BV(\Omega)\cap L^{\infty}(\Omega) ,\> \text{ with } \> \|u\|_{BV(\Omega)}\le 1 \right\}\,,
 \]
 where $\displaystyle \|u\|_{BV(\Omega)} = \int_{\partial\Omega}|u|\,d\mathcal H^{N-1} +\int_\Omega|Du|$.

 \end{Theorem}

 \begin{pf} Since we already know that
 $BV(\Omega)\cap L^\infty(\Omega)$ is a subset of $L^1(\Omega,\mu)$, all we
 have to prove is
 \begin{equation}\label{for}
 \Big|\int_\Omega u^*\,d\mu\Big|\le
 \|\mu\|_{W^{-1,\infty}(\Omega)}\Big(|Du|(\Omega)+\int_{\partial\Omega}|u|\,d\mathcal H^{N-1}\Big)\,.
 \end{equation}
 for all $u\in BV(\Omega)\cap L^\infty(\Omega)$. This inequality
 will be proved in two steps.

 Step 1: Assume first that $u\in W^{1,1}(\Omega)\cap L^\infty(\Omega)$.
 Consider the sequence $(w_n)_n$ in
 $W^{1, 1}(\Omega)\cap C(\Omega)$ of the above Lemma. Hence,
 $w_n\in L^\infty(\Omega)$ and $\|w_n\|_\infty\le\|u{\big|_{\partial\Omega}}\|_\infty$
 for all $n\in\N$. Then it yields
 \begin{multline*}
 \Big|\int_\Omega (u^*-w_n^*)\,d\mu\Big|=\big|\langle\mu,
 (u-w_n)\rangle_{W^{-1,\infty}(\Omega),W_0^{1,1}(\Omega)}\big|\le
 \|\mu\|_{W^{-1,\infty}(\Omega)}\int_\Omega|\nabla u-\nabla w_n|\,dx\\
 \\
 \le \|\mu\|_{W^{-1,\infty}(\Omega)}\Big(\int_\Omega|\nabla
 u|\,dx+\int_{\partial\Omega}|u|\,d\mathcal
 H^{N-1}+\frac1n\Big)\,.
 \end{multline*}
 It follows that
 \begin{multline}\label{uno}
 \Big|\int_\Omega u^*\,d\mu\Big|\le \Big|\int_\Omega
 (u^*-w_n^*)\,d\mu\Big|+\Big|\int_\Omega w_n^*\,d\mu\Big|\\
 \\
 \le\|\mu\|_{W^{-1,\infty}(\Omega)}\Big(\int_\Omega|\nabla
 u|\,dx+\int_{\partial\Omega}|u|\,d\mathcal
 H^{N-1}+\frac1n\Big)+\Big|\int_\Omega w_n^*\,d\mu\Big|\,.
 \end{multline}
 Since the sequence $(w_n)_n$ tends pointwise to $0$  and it is uniformly bounded
  in $L^\infty(\Omega)$, by Lebesgue's Theorem,
 $$
 \lim_{n\to\infty}\int_\Omega w_n^*\,d\mu=0\,.
 $$
 Now, taking the limit in \eqref{uno} we obtain \eqref{for}.

Step 2: In the general case, we apply Proposition \ref{M-S} and
 find a sequence
  $u_n$ in $W^{1 ,1}(\Omega)\cap C^\infty(\Omega)\cap L^\infty (\Omega)$
  such that
  \begin{equation*}
  \begin{array}{ll}
  \hbox{\rm (1) } u_n^*\to u^*\quad\hbox{ in } L^1(\Omega,\mu)\,.\\
  \\
  \hbox{\rm (2) } \int_\Omega|\nabla u_n|\,dx\to|Du|(\Omega)\,.\\
  \\
  \hbox{\rm (3) } u_n|_{\partial\Omega}=u|_{\partial\Omega}\hbox{ for all }
  n\in\N\,.\\
  \\
  \hbox{\rm (4) }
  |u_n(x)|\le\|u\|_\infty\ \ |\mu|\hbox{--a.e. for all }
  n\in\N\,.
  \end{array}
  \end{equation*}
  Then, it follows from
  $$
  \Big|\int_\Omega u_n^*\,d\mu\Big|\le
 \|\mu\|_{W^{-1,\infty}(\Omega)}\Big(\int_\Omega|\nabla u_n|\,dx+\int_{\partial\Omega}|u|\,d\mathcal
 H^{N-1}\Big)\quad\hbox{for all }n\in\N
 $$
 that \eqref{for} holds.
 \end{pf}

 \begin{Corollary}\label{clave}
 Let $\z\in\DM(\Omega)$ satisfy $\Div\z=\nu+f$ for a certain Radon measure $\nu$ and a certain $f\in L^{N,\infty}(\Omega)$. If either $\nu\ge 0$ or $\nu\le 0$,
  then $\mu=\hbox{\rm div\,}\z$ can be extended to $BV(\Omega)$ and
 $$
 \|\mu\|_{W^{-1,\infty}(\Omega)}=
 \sup\left\{\Big|\int_\Omega u^*\,d\mu\Big|\,:\, u\in BV(\Omega),
 \,|Du|(\Omega)+\int_{\partial\Omega}|u|\,d\mathcal H^{N-1}\le1\right\}\,.
 $$

  Moreover, $BV(\Omega)\hookrightarrow L^1(\Omega,\mu)$.
 \end{Corollary}

 \begin{pf}
  Consider $u\in BV(\Omega)$, denote $u_+=\max\{u,0\}$ and, for every $k>0$, apply the previous result to $T_k(u_+)$ (recall \eqref{trun}).
 Then
   \begin{multline}\label{trunc}
  \Big|\int_\Omega T_k(u_+)^*\,d\mu\Big|\le
 \|\mu\|_{W^{-1,\infty}(\Omega)}\Big(|DT_k(u_+)|(\Omega)+\int_{\partial\Omega}T_k(u_+)\,d\mathcal
 H^{N-1}\Big)\\
 \\\le
 \|\mu\|_{W^{-1,\infty}(\Omega)}\Big(|Du_+|(\Omega)+\int_{\partial\Omega}u_+\,d\mathcal
 H^{N-1}\Big)\,.
 \end{multline}
 On the other hand, observe that $u^*$ is a $\nu$--measurable
 function, so that we obtain
 $$
 \int_\Omega T_k(u_+)^*\,d\mu=\int_\Omega T_k(u_+)^*\,d\nu+\int_\Omega
 T_k(u_+(x))f(x)\,dx
 $$
 for every $k>0$. We may apply Levi's Theorem and Lebesgue's
 Theorem to deduce
 \begin{gather*}
    \lim_{k\to+\infty}\int_\Omega T_k(u_+)^*\,d\nu=\int_\Omega
    (u_+)^*\,d\nu
 \end{gather*}
 and
 \begin{gather*}
    \lim_{k\to+\infty}\int_\Omega T_k(u_+(x))f(x)\,dx=\int_\Omega
    u_+(x)f(x)\,dx\,.
\end{gather*}
  Thus,
  $$
  \lim_{k\to+\infty}\int_\Omega T_k(u_+)^*\,d\mu=\int_\Omega
    (u_+)^*\,d\mu\,.
  $$
 Now, taking the limit when $k$ goes to $\infty$ in \eqref{trunc}, it yields
 \begin{equation}\label{ec:1}
 \Big|\int_\Omega (u_+)^*\,d\mu\Big|\le
 \|\mu\|_{W^{-1,\infty}(\Omega)}\Big(|Du_+|(\Omega)+\int_{\partial\Omega}u_+\,d\mathcal
 H^{N-1}\Big)\,.
 \end{equation}

 Assume, in order to be concrete, that $\nu\ge0$.
 Since $\int_\Omega (u_+)^*\,d\mu^-=\int_\Omega
    u_+(x)f_-(x)\,dx$, we already have that $(u_+)^*$ is $\mu^-$--integrable.
    Hence, as a consequence of \eqref{ec:1}, we deduce that $(u_+)^*$ is $\mu^+$--integrable as well and then, $(u_+)^*$ $\mu$--integrable too.

 Since we may prove a similar inequality to $u_-=\max\{-u,0\}$, adding both inequalities we deduce that
 $u^*$ is $\mu$--integrable and that
 $$
 \Big|\int_\Omega u^*\,d\mu\Big|\le
 \|\mu\|_{W^{-1,\infty}(\Omega)}\Big(|Du|(\Omega)+\int_{\partial\Omega}|u|\,d\mathcal
 H^{N-1}\Big)
 $$
 holds true.
 \end{pf}

 \subsection{A Green's formula}

 Let $\z\in \DM(\Omega)$ and let $u\in BV(\Omega)$.
 Assume that $\Div\z=\nu+f$, with $\nu$ a Radon measure satisfying either $\nu\ge0$ or $\nu\le0$, and $f\in L^{N,\infty}(\Omega)$.
 In the spirit of \cite{An}, we define the following
 distribution on $\Omega$. For every $\varphi\in
C_0^\infty(\Omega)$, we write
\begin{equation}
 \label{dist1}   \langle(\z, Du),\varphi\rangle=
 -\int_\Omega u^*\,\varphi\,d\mu-\int_\Omega u\, \z\cdot\nabla
    \varphi\, dx\,,
\end{equation}
where $\mu=\Div\z$.
Note that the previous subsection implies that every term in the
above definition has sense. We next prove that this distribution
is actually a Radon measure having finite total variation.

\begin{Proposition}\label{prop}
Let $\z$ and $u$ be as above.
 The distribution $(\z, Du)$ defined previously satisfies
  \begin{equation}\label{ec:2}
  |\langle   (\z, Du), \varphi\rangle| \le \|\varphi\|_\infty \| \z
\|_{L^{\infty}(U)} \int_{U} |Du|
  \end{equation}
for all open set $U \subset \Omega$ and for all $\varphi\in
C_0^\infty(U)$.
 \end{Proposition}

  \begin{pf}
  If $U\subset \Omega$ is an open set and $\varphi\in
 C_0^\infty(U)$, then it was proved in \cite{MST2} that
  \begin{equation}\label{ec:3}
      |\langle(\z, DT_k(u)),\varphi\rangle|\le\|\varphi\|_\infty
  \|\z \|_{L^\infty(U)}\int_U|DT_k(u)|\le\|\varphi\|_\infty
  \|\z \|_{L^\infty(U)}\int_U|Du|
\end{equation}
 holds for every $k>0$.
 On the other hand,
 $$
 \langle(\z, DT_k(u)),\varphi\rangle=-\int_\Omega T_k(u)^*\varphi
 \, d\mu-\int_\Omega T_k(u)\z\cdot\nabla\varphi\, dx\,.
 $$
 We may let $k\to\infty$ in each term on the right hand side, due
 to $u^*\in L^1(\Omega, \mu)$ and $u\in L^1(\Omega)$.
 Therefore,
 $$
 \lim_{k\to\infty}\langle(\z, DT_k(u)),\varphi\rangle=\langle(\z,
 Du),\varphi\rangle\,,
 $$
 and so \eqref{ec:3} implies \eqref{ec:2}.
 \end{pf}

 \begin{Corollary}
  The distribution $(\z, Du)$ is a Radon measure. It and its total variation $\vert (\z, Du) \vert$ are absolutely
continuous with respect to the measure $\vert Du \vert$ and
$$\left\vert \int_{B}  (\z, Du) \right\vert \leq \int_{B} \vert (\z, Du) \vert \leq \Vert \z
\Vert_{L^{\infty}(U)} \int_{B} \vert Du \vert
$$
holds for all Borel sets $B$ and for all open sets $U$ such that
$B \subset U \subset \Omega$.
 \end{Corollary}

 On the other hand, for every $\z \in \mathcal{DM}^{\infty}(\Omega)$,
 a weak trace on $\partial \Omega$ of the
normal component of  $\z$ is defined in \cite{An} and denoted by
$[\z, \nu]$.

\begin{Proposition}
Let $\z$ and $u$ be as above.
 With the above definitions, the
 following Green formula holds
\begin{equation}\label{GreenI}
\int_{\Omega} u^* \, d\mu + \int_{\Omega} (\z, Du) =
\int_{\partial \Omega} [\z, \nu] u \ d\mathcal H^{N-1}\,,
\end{equation}
where $\mu=\Div\z$.
 \end{Proposition}

\begin{pf}
Applying the Green formula proved in \cite{MST2}, we obtain
  \begin{equation}\label{ec:4}
\int_{\Omega} T_k(u)^* \, d\mu + \int_{\Omega} (\z, DT_k(u)) =
\int_{\partial \Omega} [\z, \nu] T_k(u) \ d\mathcal H^{N-1}\,,
\end{equation}
for every $k>0$.  Note that the same argument appearing in the proof of the previous Proposition leads to
$$
\lim_{k\to\infty}\int_{\Omega} (\z, DT_k(u)) =\int_{\Omega} (\z,
Du)\,.
$$
We may take limits in the other terms since $u^*\in L^1(\Omega,
\mu)$ and $u\in L^1(\partial\Omega)$. Hence, letting $k$ go to
$\infty$ in \eqref{ec:4}, we get \eqref{GreenI}.
\end{pf}

\begin{Proposition}
Let $\z \in \DM(\Omega)$ with $\|\z\|_\infty \le 1$ and let $u\in BV(\Omega)$. Then $(\z,Du)=|Du|$ as measures if and only if $(\z,DT_k(u)) = |DT_k(u)|$ as measures for all $k>0$.
\label{prop1}
\end{Proposition}
\begin{pf}
We first assume $(\z,Du)=|Du|$ and so (recall \eqref{G-k-fun})
\[
    \begin{array}{rl}
    |Du| =& (\z, Du) = (\z , DT_k(u)) + (\z, DG_k(u))\\[0.2cm]
    \le & |DT_k(u)| + |DG_k(u)| = |Du|\,.\\
    \end{array}
\]
Then, the inequality becomes equality and so $(\z , DT_k(u))=|DT_k(u)|$ as measures.
\\
Conversely, we assume $(\z,DT_k(u)) = |DT_k(u)|$ for all $k>0$. For each $\varphi \in C^\infty_0 (\Omega)$, we use the same argument which appears in Proposition \ref{prop} to obtain:
\[
    \lim_{k \to \infty} \langle (\z, DT_k(u)), \varphi \rangle = \langle (\z, Du), \varphi \rangle
\]
and
\[
    \lim_{k \to \infty} \int_\Omega \varphi  \, |DT_k(u)| = \int_\Omega \varphi \, |Du|\,.
\]
So, using the hypothesis, we conclude $\langle (\z, Du), \varphi \rangle=\int_\Omega \varphi \, |Du|$ for every $\varphi \in C^\infty_0 (\Omega)$, that is, $(\z, Du)=|Du|$ as measures.
\end{pf}

\subsection{The chain rule}

We point out that there is a chain rule for $BV$--functions, the more general formula is due to L. Ambrosio and G. Dal Maso (see \cite[Theorem 3.101]{AFP}, see also \cite[Theorem 3.96]{AFP}). In our framework, it states that if $v\in BV(\Omega)$ satisfies $D^jv=0$ and $u=G(v)$, where $G$ is a Lipschitz--continuous real function, then $u\in BV(\Omega)$ and
 \[
 Du=G^\prime(v)|Dv|\,.
 \]
 We cannot directly apply this result in our context since $G^\prime$ need not be bounded. Hence, the following slight generalization is needed.

\begin{Theorem}\label{regla-cadena}
Let $v \in BV(\Omega)$ such that $D^jv=0$ and let $g$ be a continuous and unbounded real function with $g(s) > m>0$ for all $s \in \R$. We define
\[
    G(s) =\int_0^s g(\sigma)\, d\sigma\,.
\]
Assuming that $u=G(v) \in L^1(\Omega)$, it holds that $u\in BV(\Omega)$ if and only if $g(v)^*|Dv|$ is a finite measure and in that case $|Du| =g(v)^*|Dv|$ as measures.
\end{Theorem}
\begin{pf}
Let $\varphi \in C^\infty_0 (\Omega)$ with $\varphi \ge 0$.
We apply the chain rule to get the next equality:
\[
    \int_{\{ v < k\}} \varphi \, |Du| = \int_{\{ v < k\}} \varphi \, g(T_k(v))^*\,|Dv| = \int_{\{ v < k\}}\varphi \, g(v)^* \, |Dv|\,.
\]
Now, using the monotone convergence theorem, we take limits when $k \to \infty$ and it holds
\[
    \int_\Omega \varphi \, |Du| =\int_\Omega \varphi\, g(v)^*\,|Dv|\,,
\]
and if one integral is finite, the other is finite too. Finally, we generalize this equality to every $\varphi \in C^\infty_0 (\Omega)$  and the result is proved.
\end{pf}
\section{Solutions for $L^{N,\infty }$--data}

This section is devoted to solve problem
\begin{equation}
    \left \{
        \begin{array}{cl}
            \displaystyle-\Div \left( \frac{Du}{|Du|}\right) +|D u| = f(x) & \qin \Omega\,,\\
            u=0 & \qon \partial \Omega\,,
        \end{array}
    \right .
    \label{problema1}
\end{equation}
for nonnegative data $f\in L^{N,\infty}(\Omega)$. We begin by introducing the notion of solution to this problem.

\begin{Definition}\label{def1}
Let $f\in L^{N,\infty}(\Omega)$ with $f\ge0$.  We say that $u\in BV(\Omega)$ satisfying $D^ju=0$ is a \textbf{weak solution} of problem \eqref{problema1}
if there exists $\z\in\DM(\Omega)$ with $\|\z\|_{\infty} \le 1$ such that
\[
    -\Div \z + |Du| =f \,\text{ in }\, \dis (\Omega)\,,
\]
\[
    (\z, Du)=|Du| \,\text{ as measures in }\, \Omega\,,
\]
and
\[
    u\big|_{\partial \Omega} = 0\,.
\]
\end{Definition}

\begin{remark}\label{exp}\rm
We explicitly remark that any solution to problem \eqref{problema1} satisfies
\[
-\Div\big(e^{-u} \z\big) =e^{-u}f
\]
in the sense of distributions (see \cite[Remark 3.4]{MS}).
\end{remark}

\begin{Theorem}
There is a unique weak solution of problem \eqref{problema1}.
\label{teoexist}
\end{Theorem}
\begin{pf}
The proof will be divided in several steps.

\bigskip{\sl Step 1: Approximating problems.}

The function $f$ is in $ L^{N,\infty}(\Omega)$ so, there exists a sequence $\{f_n\}_{n=1}^{\infty}$ in $L^{\infty}(\Omega)$ such that $f_n$ converges to $f$ in $L^1(\Omega)$.
\\
In \cite{MS} it is proved that there exists $u_n \in
BV(\Omega)\cap L^\infty(\Omega)$, with $D^ju_n=0$ and $u_n\ge 0$, which is a solution to
problem
\begin{equation}
    \left \{
        \begin{array}{cl}
            \displaystyle-\Div \left( \frac{Du_n}{|Du_n|}\right) +|D u_n| = f_n(x) & \qin \Omega\,,\\
            u_n=0 & \qon \partial \Omega\,.
        \end{array}
    \right .
    \label{probleman}
\end{equation}
That is, there exists a vector field $\z_n$ in $\DM (\Omega)$ such that
\begin{equation}\label{condn1}
     -\Div \z_n + |Du_n| =f_n \,\text{ in }\, \dis (\Omega)\,,
\end{equation}
\begin{equation}\label{condn2}
        (\z_n, Du_n)=|Du_n| \,\text{ as measures in }\, \Omega\,,
\end{equation}
and
\begin{equation}\label{condn3}
            u_n\big|_{\partial \Omega} = 0\,.
\end{equation}
On account of Remark \ref{exp}, it also holds
\begin{equation}\label{condn4}
     -\Div (e^{-u_n}\z_n) =e^{-u_n}f_n \,\text{ in }\, \dis (\Omega)\,.
\end{equation}

\bigskip{\sl Step 2: $BV$--estimate.}

Taking the function test $\frac{T_k(u_n)}{k}$ in problem
\eqref{probleman}, we get
\[
    \frac{1}{k}\int_\Omega (\z_n,DT_k(u_n)) + \frac{1}{k}\int_\Omega T_k(u_n)^*|Du_n| = \int_\Omega f_n \frac{T_k(u_n)}{k}\, dx \le \int_\Omega f_n\, dx \le C\,,
\]
where $C$ does not depend on $n$. Since $(\z_n,Du_n) =|Du_n|$, it
follows from Proposition \ref{prop1} that $(\z_n,DT_k(u_n))=|DT_k(u_n)|$, which is nonnegative. Thus
\[
    \frac{1}{k}\int_\Omega T_k(u_n)^*|Du_n| \le C\,.
\]
Then, letting $k \to 0$ in the inequality above we arrive at
\[
    \int_\Omega |Du_n| \le C\,.
\]
Therefore, $u_n$ is bounded in $BV(\Omega)$ and, up to a
subsequence, $u_n \to u$ in $L^1(\Omega)$ and $D u_n$ converges to
$Du$ $*$--weakly as measures when $n \to \infty$.

\bigskip{\sl Step 3: Vector field.}

Now, we want to find a vector field $\z\in \DM(\Omega)$ with $\|\z\|_{\infty} \le 1$ such that
\[
    -\Div \z + |Du| \le f \,\text{ in }\, \dis (\Omega)\,.
\]
\\
The sequence $\{\z_n\}_{n=1}^{\infty}$ is bounded in
$L^\infty(\Omega; \R^N)$ then, there exists $\z\in L^\infty
(\Omega;\R^N)$ such that $\z_n \rightharpoonup \z$ $*$--weakly in
$L^\infty(\Omega;\R^N)$. In addition, since $ \|\z_n\|_\infty \le
1$ we get $\|\z\|_\infty\le 1$.
\\
Using $\varphi \in C^\infty_0 (\Omega)$ with $\varphi \ge 0$ as a function test in \eqref{probleman}, we arrive at
\[
    \int_\Omega \z_n \cdot \nabla \varphi\, dx +\int_\Omega \varphi\, |Du_n| = \int_\Omega f_n\, \varphi\, dx\,,
\]
and when we take $n \to \infty$, using \eqref{semc} it becomes
\[
    \int_\Omega \z \cdot \nabla\varphi\, dx +\int_\Omega \varphi\, |Du| \le \int_\Omega f \,\varphi\, dx\,.
\]
Therefore,
\[
    -\Div \z + |Du| \le f \qin \dis (\Omega)
\]
and $-\Div \z$ is a Radon measure. In addition, since $-\Div \z _n= f_n -|Du_n|$ holds for every $n \in \N$, the sequence $-\Div \z_n$ is bounded in
the space of measures and, due to $-\Div \z_n$ converges to $-\Div \z$, we
deduce that $-\Div \z$ is a Radon measure with finite total
variation.
\\
On the other hand, multiply \eqref{condn4} by $e^{-u_n}\varphi$, with $\varphi \in
C^\infty_0 (\Omega)$, then Green's formula provides us
\[
    \int_\Omega e^{-u_n} \z_n\cdot\nabla\varphi\, dx = \int_\Omega f_n e^{-u_n}\varphi\, dx\,,
\]
and letting $n$ go to $\infty$ we get
\[
    \int_\Omega e^{-u} \z\cdot\nabla\varphi\, dx = \int_\Omega f e^{-u}\varphi\, dx\,.
\]
Namely,
\begin{equation}\label{ecu0}
    -\Div (e^{-u}\z ) = fe^{-u}\,,\quad\hbox{in }\mathcal D'(\Omega)\,.
\end{equation}

\bigskip{\sl Step 4: $D^ju=0$.}

In this step, we are adapting an argument used in \cite{GMP}, which relies on \cite[Proposition 3.4]{AmCM} and \cite[Lemma 5.6]{C} (see also \cite[Proposition 2]{ADS}).
 A previous result is needed, namely,  inequality \eqref{desig} bellow. To prove \eqref{desig}, we begin by recalling
\[
	 -\Div (e^{-u_n} \z_n)=e^{-u_n}f_n \quad\hbox{in }\mathcal D'(\Omega)\,,
\]
since $u_n$ is the solution to problem \eqref{probleman}. Using that $u_n = G_k(u_n)+T_k(u_n)$, we can write
\[
	-\Div (e^{-u_n} \z_n)=-e^{-G_k(u_n)}\Div(e^{-T_k(u_n)}\z_n) + (e^{-u_n})^*|DG_k(u_n)|\,,
\]
and so
\begin{equation}\label{ecu1}
	\begin{array}{rcl}
	e^{-T_k(u_n)}f_n &=& -\Div (e^{-T_k(u_n)} \z_n) + (e^{-T_k(u_n)})^*|DG_k(u_n)| \\[0.2cm]
	&=& -\Div (e^{-T_k(u_n)}\z_n)+e^{-k}|DG_k(u_n)|\,.
	\end{array}
\end{equation}
Applying first the chain rule and then \cite[Proposition 2.3]{MS}, we have
\begin{multline}\label{ecu2}
	|De^{-T_k(u_n)}| = ( e^{-T_k(u_n)})^* |DT_k(u_n)|\\
	= ( e^{-T_k(u_n)})^* (\z_n,DT_k(u_n))
	= ( e^{-T_k(u_n)}\z_n,DT_k(u_n))\,.
\end{multline}
\medskip
Let $\varphi \in C^\infty_0(\Omega)$ with $\varphi \ge 0$, due to \eqref{ecu2} and \eqref{ecu1}, we get
\begin{multline*}
	\int_\Omega \varphi \, |D e^{-T_k(u_n)}| = \langle( e^{-T_k(u_n)}\z_n,DT_k(u_n)),\varphi\rangle\\
=-\int_\Omega T_k(u_n)\, \varphi\, \Div (e^{-T_k(u_n)}\z_n) - \int_\Omega T_k(u_n)\, e^{-T_k(u_n)}\,\z_n \cdot \nabla \varphi\, dx\\
	= \int_\Omega T_k(u_n)\, \varphi \, e^{-T_k(u_n)}f_n\, dx - \int_\Omega k\, e^{-k} \varphi \, |DG_k(u_n)| -\int_\Omega T_k(u_n)\, e^{-T_k(u_n)}\, \z_n \cdot \nabla \varphi\, dx\,.
\end{multline*}
That is,
\begin{multline*}
	\int_\Omega \varphi \, |De^{-T_k(u_n)}| +\frac{k}{e^k} \int_\Omega \varphi \, |DG_k(u_n)|\\
	= \int_\Omega T_k(u_n) \, \varphi \, e^{-T_k(u_n)} \, f_n \, dx- \int_\Omega T_k(u_n) \, e^{-T_k(u_n)}\,\z_n \cdot \nabla \varphi\, dx\,.
\end{multline*}
Now, we can take limits when $n$ goes to $\infty$, and applying the lower semicontinuity of the total variation, we arrive to the next inequality:
\begin{multline*}
	\int_\Omega \varphi \, |De^{-T_k(u)}| +\frac{k}{e^k} \int_\Omega \varphi \, |DG_k(u)|\\
	\le \int_\Omega T_k(u) \,\varphi \, e^{-T_k(u)} \, f\, dx - \int_\Omega T_k(u) \, e^{-T_k(u)}\, \z \cdot \nabla \varphi\, dx\,.
\end{multline*}
Finally, letting $k\to \infty$ it holds that
\[
	\int_\Omega \varphi \, |De^{-u}| \le \int_\Omega u \, \varphi \, e^{-u} \, f\, dx - \int_\Omega u \, e^{-u}\, \z \cdot \nabla \varphi\, dx
	= \langle (e^{-u}\z, Du), \varphi \rangle\,.
\]
Therefore,
\begin{equation}\label{desig}
	 |De^{-u}| \le (e^{-u}\z, Du)
\end{equation}
as measures in $\Omega$.
\\
On the other hand, we already know that
\[
	\Div (u\,e^{-u} \z) = (e^{-u} \z, Du) + u \,\Div (e^{-u}\z)\,,
\]
as measures and now we are considering the restriction on the set $J_u$. Since, by \eqref{ecu0} we have
\[
	u \,\Div (e^{-u}\z)=-u\, e^{-u}f\in L^1(\Omega)
\]
and $|J_u|=0$,
it follows that the measure $u\,\Div (e^{-u}\z)$ vanishes on $J_u$, so that
\[
	\Div (u\,e^{-u}\z)\res J_u = (e^{-u} \z, Du) \res J_u \ge |De^{-u}|\res J_u\,.
\]
Applying \cite[Lemma 2.3 and Lemma 2.4]{GMP}, the following manipulations can be performed on $J_u$:
	\begin{equation}\label{salv}
	\begin{array}{rcl}
		\Div (u\,e^{-u}\z)&=& [ue^{-u}\z, \nu_u]^+ - [ue^{-u}\z, \nu_u]^-\\[0.2cm]
		&=& u^+[e^{-u}\z,\nu_u]^+- u^-[e^{-u}\z,\nu_u]^-\,.
	\end{array}
	\end{equation}
Moreover, we also deduce that, on $J_u$,
\[
	\Div\big(e^{-u}\z\big)=[e^{-u}\z,\nu_u]^+-[e^{-u}\z,\nu_u]^-
\]
and, due to
\[
	\Div\big(e^{-u}\z\big)\in L^1(\Omega)\quad \hbox{and}\quad |J_u|=0,
\]
it follows that $[e^{-u}\z,\nu_u]^+=[e^{-u}\z,\nu_u]^-$. We will write this common value as $[e^{-u}\z,\nu_u]$. With this notation, \eqref{salv} becomes
\[
	\begin{array}{rcl}
		\Div (u\,e^{-u}\z)&=& (u^+-u^-)[e^{-u}\z,\nu_u] \\[0.2cm]
		&=& (u^+-u^-)  e^{-u^+} [\z,\nu_u]\\[0.2cm]
		&\le&(u^+-u^-)e^{-u^+}
	\end{array}
\]
Thus, we have seen that
\[
	(u^+-u^-)e^{-u^+}\mathcal H^{N-1}\res J_u\ge |De^{-u}|\res J_u =\big(e^{-u^-}-e^{-u^+}\big)\mathcal H^{N-1}\res J_u\,.
\]
Hence, for $\mathcal H^{N-1}$--almost all $x\in J_u$, we may use the Mean Value Theorem to get
\[
	(u(x)^+-u(x)^-)e^{-u(x)^+} \ge e^{-u(x)^-}-e^{-u(x)^+} = (u(x)^+-u(x)^-)e^{-w(x)}
\]
with $u(x)^-<w(x)< u(x)^+$. Therefore, it yields $u(x)^+=u(x)^-$. Since this argument holds for $\mathcal H^{N-1}$--almost every point $x\in J_u$, we get
\[
	D^ju=0\,.
\]

\bigskip{\sl Step 5: $u$ is a solution to problem \eqref{problema1}.}

To finish the proof, it remains to check that $u$ satisfies the three conditions of the definition of solution. The previous step will be essential in this checking. Indeed, it allows us to perform the following calculations:
\[
    \begin{array}{rl}
    fe^{-u} =& -\Div (e^{-u}\z ) = -(\z, D(e^{-u})^*) -(e^{-u})^*\Div \z \\[0.2cm]
    \le &   |De^{-u}|+fe^{-u} -(e^{-u})^*|Du|\\[0.2cm]
    =&fe^{-u}\,.
    \end{array}
\]
Therefore, the inequality becomes equality and so
\begin{equation}
    -\Div \z + |Du| = f \qin \dis (\Omega)\,.
    \label{distribuciones}
\end{equation}
To prove that $(\z,Du) = |Du|$ as measures in $\Omega$, we just take into account \eqref{desig}, \cite[Proposition 2.3]{MS} and the chain rule to get
\begin{equation*}
  |D(e^{-u})|\le (e^{-u}\z , Du)=(e^{-u})^*(\z, Du)\le (e^{-u})^*|Du|=|D(e^{-u})|\,,
\end{equation*}
from where the equality $(e^{-u})^*(\z, Du)= (e^{-u})^*|Du|$ as measures follows. We conclude that
$(\z,Du) = |Du|$ as measures.
\\
Now, we will prove that $u(x) = 0$ for $\mathcal H^{N-1}$--almost
all $x \in \partial \Omega$. To do that, we use the test function
$T_k(u_n)$ in problem \eqref{probleman}, so that
\[
    \int_\Omega (\z_n, DT_k(u_n)) + \int_\Omega (T_k(u_n))^*|Du_n| = \int_\Omega f \, T_k(u_n)\, dx\,.
\]
Defining the auxiliary function $J_k$ by
\[
    J_k (s) = \int_0^s T_k(\sigma)\, d\sigma =
    \left\{
        \begin{array}{lcl}
        \frac{s^2}{2}&\hbox{if}& 0\le s\le k\,,\\
        ks-\frac{k^2}{2}&\hbox{if}&k> s\,,\\
        \end{array}
    \right .
\]
we obtain
\begin{multline*}
    \int_\Omega |DT_k(u_n)| + \int_{\partial \Omega} |T_k(u_n)|\, d\mathcal H^{N-1}+ \int_\Omega |DJ_k(u_n)| +\int_{\partial \Omega} |J_k(u_n)|\, d\mathcal H^{N-1}\\
    = \int_\Omega f \,T_k(u_n)\, dx\,.
\end{multline*}
Taking into account that $J_k(u_n) \to J_k(u)$ in $L^1(\Omega)$, we let $n \to
\infty$ and applying the lower semicontinuity of functional
\eqref{semcon} we arrive at
\begin{multline*}
    \int_\Omega |DT_k(u)| + \int_{\partial \Omega} |T_k(u)|\, d\mathcal H^{N-1}+ \int_\Omega |DJ_k(u)| +\int_{\partial \Omega} |J_k(u)|\, d\mathcal H^{N-1}\\
    \le \int_\Omega f\, T_k(u)\, dx\le \int_\Omega f u\, dx\,.
\end{multline*}
Letting now $k\to \infty$ we obtain
\[
    \int_\Omega |Du| + \int_{\partial \Omega} |u|\, d\mathcal H^{N-1}+ \int_\Omega \Big|D\Big(\frac{u^2}{2}\Big)\Big| +\int_{\partial \Omega} \frac{u^2}{2}\, d\mathcal H^{N-1}\le \int_\Omega f u\, dx\,.
\]
On the other hand, Green's formula implies
\[
    \int_\Omega f u\, dx = -\int_\Omega u^*\Div \z +\int_\Omega u^*|Du| = \int_\Omega |Du| -\int_{\partial \Omega} u\, [\z,\nu]\,d\mathcal H^{N-1}+\int_\Omega u^*|Du| \,.
\]
Then
\[
    \int_{\partial \Omega} (|u| +u[\z,\nu])\,d\mathcal H^{N-1} +\int_{\partial \Omega}\frac{u^2}{2} \,d\mathcal H^{N-1} \le 0
\]
and for that, $u = 0$ in $\partial \Omega$.
\\
Now, using the same argument which is used in \cite{MS} we prove that there is a unique solution to our problem.
\end{pf}

\begin{Proposition}\label{triv}
The solution $u$ to problem \eqref{problema1} is trivial if and only if the function $f$ is such that $\|f\|_{W^{-1,\infty}(\Omega)} \le1$.
\end{Proposition}
\begin{pf}
Assume first that $\|f\|_{W^{-1,\infty}(\Omega)}\le 1 $ and let $u \in BV(\Omega)$ be the solution to problem \eqref{problema1}. Using the test function $T_k (u)$ in that problem we obtain
\begin{equation}\label{eq}
    \int_\Omega (\z,DT_k (u)) + \int_\Omega T_k (u)^*|Du| = \int_\Omega f\,T_k(u)\, dx \le \int_\Omega fu\, dx\,.
\end{equation}
Now, taking into account that $\int_\Omega T_k (u)^*|Du| \ge0$, it yields
\[
    \int_\Omega (\z,DT_k (u)) = \int_\Omega |DT_k (u)|  \le \int_\Omega fu\, dx\,.
\]
Finally,  letting $k \to \infty$ in \eqref{eq} and using H\"older and Sobolev's inequalities we arrive at
\[
    \int_\Omega |D u| +\int_\Omega u^*|Du| \le \int_\Omega fu\, dx \le \|f\|_{W^{-1,\infty}} \int_\Omega |Du|\le \int_\Omega |Du| \,.
\]
Then, $\int_\Omega u^*|Du| =0$ and thus, $u^*=0$ in $\Omega$ and we conclude $u(x)=0$ for almost every $x\in \Omega$.
\\
Now, we suppose that
\[
    \|f\|_{W^{-1,\infty}(\Omega)} = \sup \left\{ \int_\Omega \varphi f\, dx \; : \; \int_\Omega |\nabla \varphi|\, dx =1 , \; \varphi \in W^{1,1}_0 (\Omega) \right\}>1\,,
\]
that is, there exists $\psi \in W^{1,1}_0 (\Omega)$ such that
\[
    \int_\Omega |\nabla \psi|\, dx =1 \quad \text{ and } \quad \int_\Omega \psi f\, dx > 1\,.
\]
Finally, we use $\psi$ as a test function in \eqref{problema1}, so we get
\[
    \int_\Omega \psi \,|Du| = \int_\Omega \psi \, f\, dx - \int_\Omega \z \cdot \nabla \psi\, dx > \int_\Omega |\nabla \psi|\, dx  - \int_\Omega \z \cdot \nabla \psi\, dx  \ge 0\,.
\]
Therefore, $|Du| \not= 0$ and so $u\not=0$ in $\Omega$.
\end{pf}

\begin{remark}\rm
This phenomenon of trivial solutions for non--trivial data is usual in problems involving the $1$--Laplacian.
It is worth comparing the above result with \cite[Theorem 4.1]{MST1} (see also \cite[Theorem 4.2]{MST2}), where the Dirichlet problem for the equation
$-\Div\Big(\frac{Du}{|Du|}\Big)=f(x)$ is studied. Indeed, for such a problem it is seen that a datum satisfying $\|f\|_{W^{-1,\infty}(\Omega)} <1$ implies a trivial solution, while no $BV$--solution can exist for $\|f\|_{W^{-1,\infty}(\Omega)} >1$. Obviously, the most interesting case is when $\|f\|_{W^{-1,\infty}(\Omega)} =1$; then non--trivial solutions can be found for some data but the trivial solution always exists.
In our case, this dichotomy does not hold: for $\|f\|_{W^{-1,\infty}(\Omega)} =1$, only trivial solutions exist.
\end{remark}

To study the summability of the solution to problem \eqref{problema1}, we need the following technical result which will also be useful in Sections 6 and 7.

\begin{Lemma}
Let $u \in BV(\Omega)$ with $D^j u =0$ and let $\z$ be a vector field with $\|\z\|_\infty \le 1$ and $\Div \z = \mu + f$, where $\mu$ is a positive measure. If $G$ is an increasing and $C^1$ function and $\lim\limits_{s\to \infty} G(s)=\infty$, then, $(\z,Du)=|Du|$ implies $(\z, DG(u))=|DG(u)|$.
\label{lema}
\end{Lemma}
\begin{pf}
Since $(\z,Du)=|Du|$, we have $ (\z,DT_k(u))=|DT_k(u)|$ for all positive $k$. Using \cite[Proposition 2.2]{MS} we get $ (\z,DG(T_k(u)))=|DG(T_k(u))|$ for all $k>0$. Now, since $G(T_k(u))=T_{G(k)}G(u)$ and $\lim\limits_{s\to \infty} G(s) =\infty$ we apply Proposition \ref{prop1} to arrive at $ (\z,DG(u))=|DG(u)|$.
\end{pf}

\begin{Proposition} \label{prop-s}
If $u$ is the solution to problem \eqref{problema1}, then $u^n \in BV(\Omega)$ for all $n\in \N$. Consequently, $u\in L^q (\Omega)$ for all $1\le q<\infty$.
\end{Proposition}
\begin{pf}
We will prove the result by induction. If $u$ is the solution of problem \eqref{problema1}, then choosing the  solution itself as test function in problem \eqref{problema1}, we get
\[
    \int_\Omega|Du| + \int_\Omega u^*\,|Du| =\int_\Omega f\,u\, dx\,.
\]
Since the first integral is positive, we have
that $u^*\,|Du|$ is a finite measure. Thus, by Theorem \ref{regla-cadena} we know that $u^2 \in BV(\Omega)$ and $2\,u^* \,|Du|=|Du^2|$.
\\
Now, set $n \in \N$ and assume that $u^n \in BV(\Omega)$. Taking the test function $u^n$ in \eqref{problema1}, it yields
\[
    \int_\Omega (\z,Du^n) + \int_\Omega \big(u^n\big)^* |Du| = \int_\Omega fu^n\, dx\,.
\]
By Lemma \ref{lema} we have $(\z,Du^n)=|Du^n|\ge 0$, then the integral $\int_\Omega \big(u^n\big)^*|Du|$ is bounded and consequently $u^{n+1}\in BV(\Omega)$ by Theorem \eqref{regla-cadena}.
\end{pf}

\begin{remark}
If $f \in L^m(\Omega)$ for $m>N$, then the solution to problem
\eqref{problema1} belongs to $L^\infty (\Omega)$ (see \cite{MS}).
\end{remark}

\section{Radial solutions}

In this section we will show some radial solutions in $\Omega =
B_R(0)$ with $R>0$ for particular data in $L^{N, \infty}(\Omega)$.
In \cite[Section 4]{MS}, some examples of bounded solutions for
data $f \in L^q(\Omega)$, with $q>N$, can be found. In Example
\ref{boun} we show bounded solutions for $f \in
L^{N,\infty}(\Omega)\backslash L^N(\Omega)$, while in Example
\ref{ejemplo-bola} we show unbounded solutions. Therefore, unbounded solutions really occur.

Throughout this section, we will take $u(x)=h(|x|)$ with $h(r) \ge 0$, $h(R)=0$ and $h'(r)
\le 0$. To deal with the examples, we will consider two zones. If
$h'(r) < 0$, we know that $\z(x)=\frac{Du}{|Du|} = -\frac{x}{|x|}$, so that
$-\Div \z(x)=\frac{N-1}{|x|}$.
In the other case, $h'(r) = 0$ and then, the solution is constant
and we only have to determine the radial vector field $\z
(x)=\xi(|x|)\,x$, so that $\Div \z (x)=\xi'(|x|)|x|+N\xi(|x|)$. The
continuity of the vector field is always searched, otherwise it would has a jump and as a consequence, the measure $\Div\z$ would have a singular part concentrated on a surface of the form $|x|=\varrho$, and measure $|Du|$ would also have that singular part. Hence, it
would induce jumps on the solution.

\begin{Example}\label{boun}
\[
    \left \{
        \begin{array}{cl}
    \displaystyle -\Div \left(\frac{Du}{|Du|}\right) + |Du| = \frac{N-1}{|x|} + \frac{\lambda}{|x|^q} & \qin B_R(0)\,,\\
            u=0 & \qon \partial B_R(0)\,,
        \end{array}
    \right .
\]
with $0<q<1$ and $\lambda > 0$.
\end{Example}
First, we assume that $u$ is constant in a ring: $h'(r)=0$ for all $\rho_1<r<\rho_2$, and we
consider the vector field $\z(x)=x\,\xi(|x|)$. Then, denoting
$r=|x|$, the equation yields
\[
    -(r\xi'(r)+N\xi(r))= \frac{N-1}{r}+\frac{\lambda}{r^q} \,,
\]
which is equivalent to
\[
    -(r^{N}\,\xi (r))' = (N-1)\,r^{N-2}+\lambda \,r^{N-1-q}\,.
\]
Therefore, solving the equation we get the vector field
\begin{equation}\label{vect1}
    \z(x)=-x\,|x|^{-1}-\frac{\lambda}{N-q}x\,|x|^{-q}+Cx\,|x|^{-N}\,,\quad \rho_1<|x|<\rho_2\,,
\end{equation}
for some constant $C$.
We next see under what conditions we can find a value for this constant satisfying $\|\z\|_{\infty} \le 1$. To this end, we will distinguished three cases.
\begin{enumerate}
  \item Assuming that $0<\rho_1<\rho_2<R$ (and that $\z$ is continuous), if $|x|=\rho_1$, then
  \[
  -x\,|x|^{-1}=-x\,|x|^{-1}-\frac{\lambda}{N-q}x\,|x|^{-q}+Cx\,|x|^{-N}\,,
  \]
  and it implies $\frac{\lambda}{N-q}x\,|x|^{-q}=Cx\,|x|^{-N}$. Thus, we deduce that $C=\frac{\lambda}{N-q}\rho_1^{N-q}$.
  The same argument leads to $C=\frac{\lambda}{N-q}\rho_2^{N-q}$ when $|x|=\rho_2$. Therefore, $\rho_1=\rho_2$ and we have got a contradiction.
    \item If we assume $0<\rho_1<\rho_2=R$, then we may argue as above and find $C=\frac{\lambda}{N-q}\rho_1^{N-q}$. Substituting in \eqref{vect1}, we get
    \[
       \z(x)=-x\,|x|^{-1}-\frac{\lambda}{N-q}x\,|x|^{-q}+\frac{\lambda}{N-q}\rho_1^{N-q}x\,|x|^{-N}\,.
    \]
    Thus, condition $\|\z\|_\infty\le1$ yields
    \[
    \Big|1+\frac\lambda {N-q}\,|x|^{1-q}-\frac{\lambda}{N-q}\rho_1^{N-q}\,|x|^{1-N}\Big|\le 1\,.
    \]
    Nevertheless, this fact does not hold since $1+\frac\lambda {N-q}\,r^{1-q}-\frac{\lambda}{N-q}\rho_1^{N-q}\,r^{1-N}>1$ for $r>\rho_1$.
  \item If we assume $0=\rho_1<\rho_2<R$, then $\z\in L^\infty(\Omega;\R^N)$ implies $C=0$. So \eqref{vect1} becomes
  \[
  \z(x)=-x\,|x|^{-1}-\frac{\lambda}{N-q}x\,|x|^{-q}
  \]
  and it follows from $\|\z\|_\infty\le1$ that $\frac{\lambda}{N-q}x\,|x|^{-q}$ vanishes. Hence, $\lambda=0$ and a contradiction is obtained.
\end{enumerate}
In any case we get a contradiction, so that $h'(r)=0$ cannot hold on $]\rho_1,\rho_2[$. Hence, we take $\z(x)=-\frac{x}{|x|}$. Then, the equation becomes
\[
    -h'(r) = \frac{\lambda}{r^q}\,,
\]
and the solution satisfying the boundary condition is given by
\[
    u(x)= \frac{\lambda}{1-q}(R^{1-q}-|x|^{1-q})\,.
\]
\begin{remark}\rm
We may perform similar computations to those of the previous example to study problem
\[
    \left \{
        \begin{array}{cl}
    \displaystyle -\Div \left(\frac{Du}{|Du|}\right) + |Du| = \frac{N-1}{|x|} + \lambda & \qin B_R(0)\,,\\
            u=0 & \qon \partial B_R(0)\,,
        \end{array}
    \right .
\]
with $\lambda > 0$. Then the solution is given by
$u(x)=   \lambda(R-r)$,
with associated vector field
  $\z(x)=-\frac x{|x|}$.
\end{remark}

\begin{Example}\label{ejemplo-bola}
Consider $0<\rho\le R$.
\[
    \left \{
        \begin{array}{cl}
    \displaystyle -\Div \left(\frac{Du}{|Du|}\right) + |Du| = \frac{\lambda}{|x|}\car_{B_{\rho}(0)}(x) & \qin B_R(0)\,,\\
            u=0 & \qon \partial B_R(0)\,,
        \end{array}
    \right .
\]
with $\lambda > 0$.
\end{Example}

Two cases according to the value of $\lambda$ will be
distinguished:

\begin{itemize}
\item \textit{Case $0<\lambda  \le N-1$.}
\end{itemize}
Assuming $h'(r) <0$ for any $0\le r<R$, the vector field is given by $\z(x)=-\frac{x}{|x|}$ and the equation becomes
\[
    \frac{N-1}{r} - h'(r)=\frac{\lambda}{r}\car_{]0,\rho[}(r)\,.
\]
 When $\rho<R$, we have to distinguish two zones: where $\rho \le r\le R$ in which we get $h'(r)=(N-1)/r $, and where $0\le r < \rho$ in which we arrive at $h'(r)=(N-1-\lambda)/r$. Both expressions are nonnegative and so they are in contradiction with our hypothesis. We arrive at the same contradiction when $\rho=R$. Therefore, $h'(r)=0$ holds for all $0\le r<R$ and it follows $h(r)=0$ for all $0\le r<R$ due to the boundary condition. To obtain the field $\z(x)=\xi(|x|)\,x$ we have to consider the equation
\[
    -(r^N\xi(r))' =\lambda \,r^{N-2}\car_{]0,\rho[}(r)\,.
\]
If $0\le r<\rho$ we get the field $\xi(r)=-\lambda/(N-1)\,r^{-1}+Cr^{-N}$ but since we ask $\|\z\|_\infty\le1$, then $C=0$. On the other hand, if $\rho\le r < R$ we arrive at $\xi(r)=-Cr^{-N}$. In order to determine the value of $C$, we demand the continuity of $\xi$ and then the field becomes
\[
    \z(x)=\left\{
    \begin{array}{lcc}
    -\dfrac{\lambda}{N-1}\dfrac{x}{|x|}&if& 0\le r<\rho\,,\\
    - \dfrac{\lambda \rho^{N-1}}{N-1}\dfrac{x}{|x|^N} &if& \rho\le r < R\,.
    \end{array}
    \right .
\]
\begin{itemize}
\item \textit{Case $\lambda > N-1$.}
\end{itemize}
In the region $0\le r <\rho$, we may argue as in the above example and have a contradiction when $h'(r)=0$. So $h'(r)<0$ and the solution is given, up to constants, by \[u(x)=(N-1-\lambda) \log \Big(\frac{|x|}\rho\Big) \]
 with the vector field $\z(x)=-x/|x|$.
On the other hand, if $\rho < r<R$, we have a contradiction when $h'(r)<0$, wherewith the solution is $u(x)=0$ and the vector field is given by $\xi (r)=-C r^{-N}$. Since we have $\|\z\|_\infty = 1$ when $0\le r <\rho$, in order to preserve the continuity we require
\[
    1=\left| \z \left(\rho\right) \right| = C\rho^{-N} \rho\,.
\]
Therefore, the vector field becomes $\z(x)=-\rho^{N-1}\frac{x}{|x|^N}$ and the solution is given by
$$
u(x)=\left\{
\begin{array}{ll}
(N-1-\lambda)\log\big(\frac{|x|}\rho\big) &\hbox{if }0\le r\le\rho\,,\\
0 &\hbox{if } \rho<r<R\,.\\
\end{array}\right.
$$

\begin{remark}\label{nota}\rm
An important particular case of the previous example is the
problem
\begin{equation}\label{ejm-imp}
    \left \{
        \begin{array}{cl}
    \displaystyle -\Div \left(\frac{Du}{|Du|}\right) + |Du| = \lambda \frac{1}{|x|} & \qin B_R(0)\,,\\
            u=0 & \qon \partial B_R(0)\,,
        \end{array}
    \right .
\end{equation}
with $\lambda > 0$. We have seen that the solution is given by
$$
u(x)=\left\{
\begin{array}{ll}
0 &\hbox{when }0<\lambda  \le N-1\,,\\
(N-1-\lambda)\log\big(\frac{|x|}R\big) &\hbox{when }\lambda >N-1\,.
\end{array}\right.
$$

Problem \eqref{ejm-imp} can be seen as the limit case of
problems with a Hardy--type potential, namely,
$$    \left \{
        \begin{array}{cl}
    \displaystyle -\Div \left(|\nabla u|^{p-2}\nabla u\right) + |\nabla u|^p = \lambda\frac{u^{p-1}}{|x|^p} & \qin B_R(0)\,,\\
            u=0 & \qon \partial B_R(0)\,,
        \end{array}
    \right .
    $$
    Problems with Hardy--type potential received much attention in recent years. We point out that in
    \cite{APP} has been studied problem \eqref{ejm-imp} with $p=2$
    showing the regularizing effect produced by the gradient term as absorption.
\end{remark}

\section{Changing the unknown: More general gradient terms}

From now on, we will generalize problem \eqref{problema1} adding a continuous function $g : [0,\infty[ \to \R$  in the gradient term:
\begin{equation}
    \left \{
        \begin{array}{cl}
            \displaystyle-\Div \left( \frac{Dv}{|Dv|}\right) +g(v)\,|D v| = f(x) & \qin \Omega\,,\\
            v=0 & \qon \partial \Omega\,.
        \end{array}
    \right .
    \label{problemag}
\end{equation}
In this section, this problem will be studied for a function $g$ that will result in standard cases.

The existence and uniqueness of solutions to problem \eqref{problemag} depend on the properties of the function $g$, and the definition of solution to a problem may depend of the case we are studying. In any case, we have to give a sense to $g(v)|Dv|$, since the meaning of that term depends on the representative of $g(v)$ we are actually considering.
First of all, we will assume that a solution satisfies $D^jv=0$ and then we will take $g(v)$ as the precise representative $g(v)^*=g(v^*)$, which is integrable with respect to the measure $|Dv|$.

\subsection{Bounded $g$}
In this subsection, let $g$ be a continuous and bounded function such that there exists $m>0$ with $g(s)\ge m$ for all $s\ge0$. We define the function
\[
    G(s)=\int_0^s g(\sigma)\,d\sigma\,.
\]
With this notation, the term $g(v)|Dv|$ in the equation means $|DG(v)|$.
\begin{Definition}\label{def-g}
We say that a function $v$ is a \textbf{weak solution} to problem \eqref{problemag} with  $g$  defined as above, if $v \in BV (\Omega)$ with $D^jv=0$ and there exists a field $\z\in \DM (\Omega)$ with $\|\z\|_\infty \le 1$ such that
\[
    -\Div\z + g(v)^*|Dv| = f \qin \dis (\Omega)\,,
\]
\[
    (\z, Dv)=|Dv| \quad \text{ as measures in }\quad \Omega\,,
\]
and
\[
    v\big|_{\partial \Omega} = 0\,.
\]
\end{Definition}

\begin{Theorem}
Let $u$ be the solution to problem \eqref{problema1}. Assume that $g$ is a continuous real function such that $0<m\le g(s)$ for all $s\ge0$ and let $u=G(v)$. Then, $v$ is a solution to problem \eqref{problemag}.
\label{teosolg}
\end{Theorem}
\begin{pf}
Since the function $u$ is the solution of problem \eqref{problema1}, there exists a vector field $\z\in \DM (\Omega)$ such that
\begin{equation}
    -\Div\z + |Du| = f \qin \dis (\Omega)\,,
    \label{cond1}
\end{equation}
\begin{equation*}
    (\z, Du)=|Du| \quad \text{ as measures in }\quad \Omega\,,
    \label{cond2}
\end{equation*}
and
\[
    u\big|_{\partial \Omega} = 0\,.
\]
By the properties of $g$,  the function $G$ is increasing and the derivative of $G^{-1}$ is bounded. Then, we apply the chain rule to get $v=G^{-1}(u)\in BV(\Omega)$. We also deduce $D^jv=0$ and
\[
    v\big|_{\partial \Omega} = G^{-1}(u)\big|_{\partial \Omega} = 0\,.
\]
Moreover, it holds by Lemma \ref{lema}:
\[
    (\z,Dv)=|Dv| \quad \text{ as measures in }\quad \Omega\,.
\]
Finally, making the substitution $u=G(v)$ in \eqref{cond1} and applying the chain rule we get
\[
    -\Div\z + g(v)^*|Dv| = f \qin \dis (\Omega)\,.
\]
\end{pf}

\begin{Corollary}
If $v$ is a solution to problem \eqref{problemag} with $g$ continuous, bounded and such that $g(s)\ge m>0$ for all $s\ge0$, then, $u=G(v)$ is the solution to  problem \eqref{problema1}.
\label{corsolg}
\end{Corollary}
\begin{pf}
Applying the same argument which is used in Theorem \ref{teosolg} and keeping it in mind that $g$ is bounded and $G$ is increasing, the result is proved.
\end{pf}

\begin{Theorem}
There exists a unique solution to problem \eqref{problemag} with $g$ continuous, bounded and such that $g(s)\ge m>0$ for all $s\ge0$.
\label{teounicidad}
\end{Theorem}
\begin{pf}
Assuming there are two solutions $v_1$ and $v_2$ of problem \eqref{problemag}, by the Corollary \ref{corsolg}, $G(v_1)$ and $G(v_2)$ are solutions to problem \eqref{problema1}. Thus, $G(v_1)=G(v_2)$ and since $G$ is injective we get $v_1=v_2$.
\end{pf}

\subsection{Unbounded $g$}

In this subsection we will prove an existence and uniqueness result to problem \eqref{problemag} assuming $g(s)\ge m>0$ be an unbounded function.
\begin{Theorem}\label{teoalpha}
There is a unique solution to problem \eqref{problemag} with $g$ continuous and such that $g(s)\ge m>0$ for all $s\ge0$.
\end{Theorem}
\begin{pf}
First of all, we consider the approximate problem
\begin{equation}
    \left\{
    \begin{array}{cl}
        \displaystyle -\Div \left( \frac{Dv_k}{|Dv_k|}\right) +T_k(g(v_k))|D v_k| = f(x) & \qin \Omega\,,\\
        v_k=0 & \qon \partial \Omega\,.
    \end{array}
    \right .
    \label{problalphak}
\end{equation}
By Theorem \ref{teounicidad}, it has a unique solution. Then, there exists $v_k \in BV(\Omega)$ with $D^jv_k =0$ and also a vector field $\z_k \in \DM(\Omega)$ such that $\|\z_k\|_\infty\le 1$ and
\[
    -\Div \z_k +T_k(g(v_k))^*|Dv_k| = f \qin \dis (\Omega)\,,
\]
\[
    (\z_k,Dv_k)=|Dv_k| \quad \text{as measures}\,,
\]
and
\[
    v_k\big| =0 \quad \mathcal H^{N-1}\text{--a.e. in } \partial \Omega\,.
\]
First, we take the test function $\frac{T_h (v_k)}{h}$ in problem \eqref{problalphak} and we get
\begin{equation*}
    \frac{1}{h}\int_\Omega (\z_k, D T_h(v_k)) + \int_\Omega T_k(g(v_k))^* \frac{T_h(v_k)^*}{h} |Dv_k| = \int_\Omega f \,\frac{T_h(v_k)}{h}\, dx \le \int_ \Omega f\, dx\,.
\end{equation*}
Keeping in mind that the first integral is positive (by Lemma \ref{lema}), we can take limits in the second integral when $h \to 0$ and so we obtain
 \begin{equation}\label{eq1}
    \int_\Omega T_k(g(v_k))^* |Dv_k|  \le \int_ \Omega f\, dx\,.
\end{equation}
 Since $T_k(g(v_k))$ is bigger than $m$,  it yields
\[
    m\,\int_\Omega |Dv_k| \le \int_\Omega f\, dx\,.
\]
Therefore, $v_k$ is bounded in $BV(\Omega)$ and there exists $v\in BV(\Omega)$ such that, up to subsequences, $v_k \to v$ in $L^1(\Omega)$ and a.e.. Moreover, $Dv_k \to Dv$ $*$--weak as measures when $k \to \infty$.
\\
To prove $D^j v =0$ we use the same argument which appears in Theorem \ref{teoexist}, so we get $D^j G(v) =0$ and then we deduce that $D^j v =0$. On the other hand, we define the function
\[
    F_k(s):=\int_0^s T_k (g(\sigma))\, d\sigma\,.
\]
Using \eqref{eq1} and the chain rule we have the next inequality:
\[
    \int_\Omega |DF_k(v_k)| \le \int_\Omega f\, dx\,.
\]
which implies that the sequence $F_k(v_k)$  is bounded in $BV(\Omega) $ and converges in $L^1(\Omega)$ to $G(v)$. Now, denoting $u_k=F_k(v_k)$ and $u=G(v)$ we get that $u_k$ converges to $u$ in $L^1(\Omega)$ and
\[
    \int_\Omega|Du_k| \le \int_\Omega f\, dx \,.
\]
Therefore, it is true that $u\in BV(\Omega)$. Moreover,  keeping in mind Theorem \ref{regla-cadena}, we get $|Du|=g(v)^*|Dv|$ as well.
\\
By Corollary \ref{corsolg}, $u_k$ is the solution to
\[
    \left\{
    \begin{array}{cl}
        \displaystyle -\Div \left( \frac{Du_k}{|Du_k|}\right) +|D u_k| = f(x) & \qin \Omega\,,\\
        u_k=0 & \qon \partial \Omega\,.
    \end{array}
    \right .
\]
The same argument used in the proof of Theorem \ref{teoexist} works for determining that $u$ is the solution to
\[
    \left\{
    \begin{array}{cl}
        \displaystyle -\Div \left( \frac{Du}{|Du|}\right) +|D u| = f(x) & \qin \Omega\,,\\
        u=0 & \qon \partial \Omega\,.
    \end{array}
    \right .
\]
Finally, since $g(s)\ge m>0$ for all $s\ge0$ and applying Theorem \ref{teosolg}, we deduce that $v$ is the solution to problem \eqref{problemag}.
\end{pf}

\begin{Proposition} \label{propq}
The solution $v$ to problem \eqref{problemag} satisfies $v\in L^q (\Omega)$ for all $1\le q<\infty$.
\end{Proposition}
\begin{pf}
The proof follows the argument of the proof of Proposition \ref{prop-s}, on account of $g(s)\ge m>0$ for all $s\ge0$.
\end{pf}

\section{A non standard case: $g$ touches the axis}

In this section we assume that $g$ is a continuous, bounded and non integrable function with $g(s) >0$ for almost every $s\ge0$. In this case, $G$ is increasing but $(G^{-1})'$ may be unbounded.

First, we analyze the case when there exist $m, \, \sigma>0$ such
that $g(s)\ge m>0$ for all $s\ge \sigma $. Observe that this
condition resembles Condition (1.7) in \cite{ABPP}.

\begin{Theorem}\label{teo-finitos-ceros}
Let $g$ be as above. Then, there exists a solution to problem \eqref{problemag}.
\end{Theorem}
\begin{pf}
Let $v_n$ be the solution to the approximating problem
\[
    \left\{
    \begin{array}{cl}
        \displaystyle -\Div \left( \frac{Dv_n}{|Dv_n|}\right) +\left(g(v_n)+\frac{1}{n}\right)|D v_n| = f & \qin \Omega\,,\\
        v_n=0 & \qon \partial \Omega\,,
    \end{array}
    \right .
\]
with the associated vector field $\z_n$. Using the test function
$\frac{T_k(v_n-T_{\sigma}(v_n))}{k}$ in that problem we get
\[
\int_{\{ v_n>\sigma \}} g(v_n)^* \frac{T_k(v_n-T_{\sigma}(v_n))^*}{k}
|Dv_n| \le \int_{\{v_n>\sigma\}} f\, dx\,;
\]
and taking limits when $k\to 0^+$ it yields
\[
    \int_{\{ v_n>\sigma\}} g(v_n)^* |Dv_n| \le \int_{\{v_n>\sigma\}} f\, dx\,.
\]
Since there exist $m>0$ such that $g(s)\ge m$ for all $s\ge \sigma$, then, the previous inequality becomes:
\begin{equation}\label{parte1}
     \int_{\{v_n>\sigma\}} |Dv_n| \le \frac{1}{m} \int_\Omega f\, dx\,.
\end{equation}
Now, we use the test function $T_{\sigma}(v_n)$ in the same problem, so we get
\begin{equation}\label{parte2}
    \int_{\{v_n \le \sigma\}} |Dv_n| \le \int_\Omega f \, T_{\sigma}(u_n)\, dx \le  \sigma \int_\Omega f\, dx\,.
\end{equation}
Finally, with \eqref{parte1} and \eqref{parte2} we have
\[
    \int_\Omega |Dv_n| \le \left( \sigma +\frac{1}{m}\right) \int_\Omega f\, dx \quad \mbox{ for all } n\in \N\,,
\]
that is, the sequence $(v_n)_n$ is bounded in $BV(\Omega)$ and this implies that, up to subsequences, there exists $v\in BV(\Omega)$ with $v_n \to v$ in $L^1(\Omega)$ and a.e. as well as $Dv_n \to Dv$ $*$--weak in the sense of measures. We conclude the proof using  arguments of Theorem \ref{teoexist}.
\end{pf}

For a general function $g$ we have to change the definition of solution. We will show in Example \ref{toca-infin} that Definition \ref{def-g} does not really work.
\begin{Definition}\label{nueva}
Let $g$ be a continuous, bounded and non integrable function with $g(s)>0$ for almost every $s\ge0$. We say that a function $v$ is a \textbf{weak solution} to problem \eqref{problemag} if $v(x)<\infty$ a.e. in $\Omega$, $G(v) \in BV (\Omega)$ with $D^jG(v)=0$ and there exists a field $\z\in \DM (\Omega)$ with $\|\z\|_\infty \le 1$ such that
\[
    -\Div\z + g(v)^*|Dv| = f \qin \dis (\Omega)\,,
\]
\[
    (\z, DG(v))=|DG(v)| \quad \text{ as measures in }\quad \Omega\,,
\]
and
\[
    v\big|_{\partial \Omega} = 0\,,
\]
where the function $G$ is defined by
\[
    G(s)=\int_0^s g(\sigma) \, d\sigma\,.
\]
\end{Definition}

\begin{Theorem}\label{teoproblemag}
Assume that the function $g$ is continuous, bounded and non integrable with $g(s) >0$ for almost every $s\in \R$. Then, there exists a unique solution to problem \eqref{problemag} in the sense of Definition \ref{nueva}.
\end{Theorem}
\begin{pf}
The approximating problem
\begin{equation}
    \left \{
        \begin{array}{cl}
            \displaystyle-\Div \left( \frac{Dv_n}{|Dv_n|}\right) +\left(g(v_n)+\frac{1}{n}\right)|D v_n| = f(x) & \qin \Omega\,,\\
            v_n=0 & \qon \partial \Omega\,,
        \end{array}
    \right .
    \label{problemagn}
\end{equation}
has a unique solution for every $n\in \N$ because of Theorem \ref{teounicidad}. That is, there exists a vector field $\z_n\in \DM (\Omega)$ with $\|\z_n\|_\infty\le 1$ and a function $v_n \in BV(\Omega)$ with $D^j v_n = 0$ and such that
\begin{equation}
    -\Div\z_n + \left( g(v_n)+\frac{1}{n}\right)^*|Dv_n| = f \qin \dis (\Omega)\,,
    \label{condgn}
\end{equation}
\[
    (\z_n, DG_n(v_n))=|DG_n(v_n)| \quad \text{ as measures in }\quad \Omega\,,
\]
and
\[
    v_n\big|_{\partial \Omega} = 0\,,
\]
where we denote
\[
    G_n(s)=\int_0^s \left( g(\sigma) +\frac{1}{n}\right) \, d\sigma\,.
\]
We will show that the limit of the sequence $(v_n)_n$ is the solution to problem \eqref{problemag}. First of all, we take the test function $\frac{T_k(v_n)}{k}$ in problem \eqref{problemagn} and we arrive at
\[
    \frac{1}{k}\int_\Omega T_k(v_n)^*|DG_n(v_n)| \le \int_\Omega f\, dx
\]
for every $k$. Now, letting $k \to 0$ and using Fatou's Theorem we get
\[
     \int_{\{v_n \not =0 \}} |DG_n(v_n)|\le \int_\Omega f\, dx\,.
\]
In addition, since $D^jv_n = 0$ it follows that $Dv_n=0$ almost everywhere in $\{ v_n=0\}$. Thus,
\[
     \int_\Omega |DG_n(v_n)|\le \int_\Omega f\, dx\,,
\]
and so $G_n(v_n)$ is bounded in $BV(\Omega)$. This implies that, up to subsequences, there exist $w$ such that $G_n(v_n) \to w$ in $L^1(\Omega)$ and a.e., and also $DG_n(v_n) \to Dw$ $*$--weak in the sense of measures.
We denote $v=G^{-1}(w)$, which is finite a.e..
\\
In what follows, we apply the same argument used in Theorem \ref{teoexist} with minor modifications, hence we just sketch it. We get $\z_n \rightharpoonup \z$ $*$--weakly in $L^\infty (\Omega)$ with $\|\z\|_\infty \le 1$ and $-\Div \z$ is a Radon measure with finite total variation. Moreover, using the test function $e^{-G_n(v_n)}\varphi$ with $\varphi \in C^\infty_0 (\Omega)$ in problem \eqref{problemagn} and letting $n$ go to $\infty$, it leads $-\Div (e^{-G(v)}\z) = e^{-G(v)}f$ in the sense of distributions. The next step is to show, with the same argument used in Theorem \ref{teoexist}, that $D^jG(v)=0$ and deduce $D^jv=0$.
 Then is easy to obtain
\[
    -\Div \z + |DG(v)| = f \qin \dis (\Omega)
\]
in the sense of distributions and
\[
(\z,DG(v))=|DG(v)|
\] as measures. Moreover, we take $T_k(G_n(v_n))$ in \eqref{problemagn} to arrive at $G(v)\big|_{\partial \Omega}=0$ and then, we also get
\[
    v\big|_{\partial \Omega} = 0\,.
\]

The uniqueness can be proved as in \cite{MS}.
\end{pf}

To remark the necessity to have a new definition to the concept of solution, we show in the next example that the solution to \eqref{problemag} when $g$ is such that $\lim\limits_{s\to \infty}g(s)=0$ is not in $BV(\Omega)$.

\begin{Example}\label{toca-infin}
The solution to problem
\begin{equation}
    \left \{
        \begin{array}{cl}
            \displaystyle-\Div \left( \frac{Dv}{|Dv|}\right) +\frac{1}{1+v}|D v| = \frac{\lambda}{|x|} & \qin \Omega\,,\\
            v=0 & \qon \partial \Omega\,,
        \end{array}
    \right .
    \label{probl-toca-infin}
\end{equation}
is not in $BV(\Omega)$ for $\lambda$ big enough.
\end{Example}
First, we will solve the related problem
\begin{equation}
    \left \{
        \begin{array}{cl}
            \displaystyle-\Div \left( \frac{Du}{|Du|}\right) +|D u| = \frac{\lambda}{|x|} & \qin \Omega\,,\\
            u=0 & \qon \partial \Omega\,,
        \end{array}
    \right .
    \label{probl-relacionado}
\end{equation}
and then, using the  inverse function of
\[
    G(s)=\int_0^s \frac{1}{1+\sigma}\, d\sigma = \log (1+s)
\]
we will get the solution $v$.
\\
Due to Example \ref{ejemplo-bola} we know that, for $\lambda>N-1$, the solution to problem \eqref{probl-relacionado} is given by  $u(x)=(N-1-\lambda)\log (|x|/R)$ with  the associated field $\z(x)=-x/|x|$. Moreover, the inverse of function $G$ is given by $G^{-1}(s)= e^s-1$. Therefore, the solution to \eqref{probl-toca-infin} is given by
\[
    v(x)=G^{-1}(u(x))=\left(\frac{|x|}{R}\right)^{N-1-\lambda}-1
\]
when $\lambda >N-1$. Nevertheless, $v$ is not in $BV(\Omega)$ when $N<\lambda/2 +1$ because in that case, $\displaystyle |Du| = \frac{\lambda-N+1}{R^{N-1-\lambda}}|x|^{N-2-\lambda}$ is not integrable.

\section{Odd cases}

In this last section we will show some cases where the properties of the function $g$ does not provide uniqueness, existence or regularity of solutions to problem \eqref{problemag}.

\subsection{First case}

First of all, we suppose the function $g$ is integrable. With that condition about $g$, it is the function $f$ who determines the existence or absence of solution.

\begin{Theorem}\label{raro1}
Let $f\in L^{N, \infty}(\Omega)$ with $f\ge0$ and we consider problem \eqref{problemag} with $g\in L^1([0,\infty[)$. Then,
\begin{itemize}
\item[(i)] if $\|f\|_{W^{1,-\infty}(\Omega)} \le 1$, the trivial solution holds;
\item[(ii)] if $\|f\|_{W^{1,-\infty}(\Omega)} > e^{G(\infty)}$, does not exist any solution;
\end{itemize}
with $G(\infty)=\sup\, \{ G(t) \, :\, s\in ]0, \infty [\} $.
\end{Theorem}
\begin{pf}
The first point is deduced following the proof of Proposition \ref{triv}.
\\
On the other hand, let $\varphi \in W^{1,1}_0(\Omega)$, we use $-\Div (e^{-G(v)}\z) = e^{-G(v)}f$ to get
\[
    e^{-G(\infty)} \int_\Omega f\, |\varphi|\, dx \le \int_\Omega e^{-G(u)}f\,|\varphi|\, dx = \int_\Omega e^{-G(u)} \z \cdot \nabla |\varphi| \, dx \le \int_\Omega |\nabla  \varphi|\, dx\,.
\]
Then, if $\|f\|_{W^{-1,\infty}(\Omega)} > e^{G(\infty)}$, cannot exist any solution to problem \eqref{problemag}.
\end{pf}

\begin{remark}\rm
Since we have shown in \eqref{rel-norm} that
\[
	\|f\|_{W^{-1,\infty}(\Omega)} \le S_N \|f\|_{L^{N,\infty}(\Omega)}\,,
\]
Theorem \ref{raro1} implies the following fact:
\begin{itemize}
\item[(i)] If $\|f\|_{L^{N,\infty}(\Omega)} \le S_N^{-1}$, the trivial solution holds.
\end{itemize}
\end{remark}

\begin{remark}\rm
One may wonder what happens when $1<\|f\|_{W^{-1,\infty}(\Omega)}\le e^{G(\infty)}$. Consider the approximate solutions
$v_n$ to problem \eqref{problemagn} and let $w$ satisfy $G(v_n) \to w$. Then $w \in [0, G(\infty)]$. In particular, if $w\in [0, G(\infty)[$, the function $v = G^{-1}(w)$ is finite a.e. in $\Omega$ and is the solution to problem \eqref{problemag}. However, $w$ can be equal to $G(\infty)$ in a set of positive measure and so $v$ is infinite in the same set. We conclude that $v$, in this case, is not solution.
\end{remark}

\begin{Example}\label{no-sol-radial}
Problem
\begin{equation}
    \left \{
        \begin{array}{cl}
            \displaystyle-\Div \left( \frac{Dv}{|Dv|}\right) +\frac{1}{1+v^2}\,|D v| = \frac{N-1}{|x|}+\lambda & \qin B_R(0)\,,\\
            v=0 & \qon \partial B_R(0)\,,
        \end{array}
    \right .
    \label{nosolradial}
\end{equation}
has not radial solutions when $\lambda$ is large enough.
\end{Example}
Assuming there exists a radial solution $u(x)=h(|x|)$  with $h:[0,R] \to \R$ is such that $h(r)\ge 0$, $h(R)=0$ and $h'(r)\le 0$, we will get a contradiction.
First, we suppose that  $h'(r)=0$ for $\rho_1<r<\rho_2$ and, reasoning as in Example \ref{boun}, we get a contradiction. Therefore, we only can have $h'(r)<0$ for all $0\le r<R$. In this case, we know that the vector field is given by $\z(x)=-x/|x|$ and the equation becomes
\[
    -g(h(r))h'(r)=\lambda\,,
\]
which is equivalent to $(G(h(r))' = -\lambda$. Then, the solution is given by $G(h(r))=\lambda (R-r)$.
\\
On the other hand, we know that $G(s)\in [0,\frac{\pi}{2}[$ because
\[
    G(s) = \int_0^s g(\sigma)\,d\sigma =\int_0^s \frac{1}{1+\sigma^2}\,d\sigma = \arctan s\,.
\]
Thus, we have a radial solution if $\displaystyle\lambda < \frac{\pi}{2R}$. When $\displaystyle\lambda = \frac{\pi}{2R}$, we also obtain a radial solution, which is given by
\[
u(x)=\tan\big(\lambda(R-r)\big)\,.
\]

\subsection{Second case}

Now, we will take the function $g :[0,\infty[ \to \R$ such that
$g(s)=0$ when $s \in [0,\ell]$ and $g(s) > 0$ for all $s>\ell$. We assume $g\not\in L^1([0,\infty[)$ as well.

\begin{remark}\label{no-unicidad1}\rm
With $g$ defined as above, there is not uniqueness of solutions.

On the one hand, if $\|f\|_{L^{N,\infty}(\Omega)} \le S_N^{-1}$ and $u \in BV(\Omega)$ satisfies $u\big|_{\partial \Omega} =0$, then the function $T_\ell (u)$ is a solution to problem \eqref{problemag}. Thus, there is not uniqueness in any way.

On the other hand, if $\|f\|_{L^{N,\infty}(\Omega)} >S_N^{-1}$ we define
\[
    h(s)=g(s+\ell)
\]
and let $w$ be a solution to problem
\begin{equation}
    \left \{
        \begin{array}{cl}
            \displaystyle-\Div \left( \frac{Dw}{|Dw|}\right) +h(w)\,|D w| = f & \qin \Omega\,,\\
            w=0 & \qon \partial \Omega\,,
        \end{array}
    \right .
    \label{problemah}
\end{equation}
with associated field $\z$. Therefore, $v(x)=w(x)+\ell$ is a solution to problem \eqref{problemag} with the same vector field $\z$.


Moreover, let $\psi : [0,\ell+1] \to [\ell,\ell+1]$ be  an increasing and bijective $C^1$--function such that $\psi'(\ell+1)=1$. Then we consider
\[
    h(s)=
    \left\{
    \begin{array}{lcc}
        \psi'(s)g(\psi (s)) &if& 0 \le s\le \ell+1\,,\\
        g(s) &if& \ell+1<s\,,\\
    \end{array}
    \right .
\]
and let $w$ be a solution to problem \eqref{problemah} with $h$ defined as above. Therefore, the function
\[
    v(x)=
    \left\{
    \begin{array}{lcc}
        \psi (w(x)) &if& 0\le w(x)\le \ell+1\,,\\
        w(x) &if& \ell+1<w(x)\,,\\
    \end{array}
    \right .
\]
is a solution to \eqref{problemag}, as we can see as follows.
It is straightforward that the equation holds in $\dis(\Omega)$ and $v\big|_{\partial \Omega} =0$. We only have to see that $(\z,DG(v))=|DG(v)|$ as measures in $\Omega$.
If $0\le s\le \ell+1$ we get
\[
    H(s)=\int_0^s h(\sigma) \,d\sigma = \int_0^s \psi'(\sigma)g(\psi (\sigma))\,d\sigma = \int_0^{\psi(s)} g(\sigma) \, d\sigma = G(\psi (s))\,,
\]
\[
    H(\ell+1)=G(\psi (\ell+1)) = G(\ell+1)\,,
\]
and for $s>\ell+1$ we have
\[
    H(s) = H(\ell+1)+\int_{\ell+1}^s h(\sigma)\,d\sigma = G(\ell+1)+\int_{\ell+1}^s g(\sigma)\,d\sigma = G(s)\,.
\]
Therefore, $DG(v(x))=DH(w(x))$ and we conclude $(\z,DG(v))=|DG(v)|$ as measures in $\Omega$.
\end{remark}

\begin{Example} \label{ejemplo-no-frontera}
The solution to problem
\begin{equation}
    \left \{
        \begin{array}{cl}
            \displaystyle-\Div \left( \frac{Du}{|Du|}\right) + g(u) |Du| = \frac{N}{|x|} & \qin \Omega\,,\\
            u=0 & \qon \partial \Omega\,,
        \end{array}
    \right .
\label{ejemplog}
\end{equation}
with
\[
    g(s)=
    \left\{
        \begin{array}{lcc}
        0 &if& s\le a\,,\\
        s-a &if& a<s\,,
        \end{array}
    \right .
\]
for $a>0$ does not vanish on $\partial \Omega$.
\end{Example}
We define
\[
    G(s) = \int_0^s g(\sigma)\,d\sigma = \left\{
    \begin{array}{lcl}
    0 & \mbox{ if } & 0\le s\le a\,,\\
    \displaystyle \frac{a}{2} + \frac{s^2}{2}-a\,s & \mbox{ if } &  a<s\,.
	\end{array}
	\right .
\]
It is easy to prove that
\[
    u(x) = h(|x|)=h(r)=G^{-1} \left(-\log \left( \frac{r}{R}\right)\right)
\]
with $\z = \frac{x}{|x|}$ is such that $(\z,Du) = |Du|$ as measures in $\Omega$ and $-\Div \z +g(u)^*|Du| = \frac{N}{r}$ in $\dis (\Omega)$. However,
\[
    h(R)=G^{-1} (0) = 1\,.
\]
Although the boundary condition is not true, the solution achieves the boundary weakly (see \cite{ABCM}), that is
\[
    [\z,\nu] = -\frac{x}{|x|}\frac{x}{|x|} = -1 = -\sg (u)\,.
\]

\subsection{Third case}
Finally, let $0<a<b$, we will take $g$ a function with $g(s)=0$ when $s \in [a, b]$ and $g(s) > 0$ for all $s<a$ and $s>b$. Moreover we assume that $g\not\in L^1([0,\infty[)$.

\begin{remark}\label{no-unicidad2}\rm
We will use a similar argument to the previous one to show that there is not uniqueness of solution to problem \eqref{problemag} with function $g$ defined as above.

Let $\psi : [0,b] \to [0,a]$ be an increasing and bijective $C^1$--function. Now, we define
\[
    h(s)=
    \left\{
    \begin{array}{lcc}
        \psi'(s)g(\psi (s)) &if& 0 \le s\le b\,,\\
        g(s) &if& b<s\,.\\
    \end{array}
    \right .
\]
If $w$ is a solution to problem \eqref{problemah}, then, we have that
\[
    v(x)=
    \left\{
    \begin{array}{lcc}
        \psi(w (x)) &if& 0\le w(x)\le b\,,\\
        w(x) &if& b<w(x)\,,\\
    \end{array}
    \right .
\]
 is a solution to the original problem \eqref{problemag} because the equation holds in $\dis(\Omega)$ and also $w\big|_{\partial \Omega} = 0$. In addition, for $0\le s\le b$ we have
\[
   H(s)=\int_0^s h(\sigma) \,d\sigma = \int_0^s \psi'(\sigma)g(\psi (\sigma))\,d\sigma = \int_0^{\psi(s)} g(\sigma) \, d\sigma = G(\psi (s))\,,
\]
\[
   H(b)=G(\psi (b))= G(a) = G(b)
\]
and for $s>b$ we get
\[
   H(s) = H(b)+\int_b^s h(\sigma)\,d\sigma = G(b)+\int_b^s g(\sigma) \, d\sigma = G(s)\,.
\]
Therefore, we have proved the remaining condition: $(\z, DG(v))=|DG(v)|$ as measures in $\Omega$.
\end{remark}

\begin{Example}
\label{ejemplo-disc}
Problem
\begin{equation}
    \left \{
        \begin{array}{cl}
            \displaystyle-\Div \left( \frac{Du}{|Du|}\right) + g(u) |Du| = \frac{N}{|x|} & \qin \Omega\,,\\
            u=0 & \qon \partial \Omega\,,
        \end{array}
    \right .
\label{ejemplog}
\end{equation}
with
\[
    g(s)=
    \left\{
        \begin{array}{lcc}
        a-s&if&s<a\,,\\
        0 &if&a\le s\le b\,,\\
        s-b &if& b<s\,,
        \end{array}
    \right .
\]
where $0<a<b$, has a discontinuous solution.
\end{Example}
We define
\[
    G(s) = \int_0^s g(\sigma)\,d\sigma =\left\{
    \begin{array}{lcl}
    \displaystyle \frac{-s^2}{2}+a\,s  & \mbox{ if } & 0\le s\le a\,,\\[0.2cm]
    \displaystyle \frac{a^2}{2} & \mbox{ if } & a\le s\le b\,,\\[0.2cm]
    \displaystyle \frac{a^2+b^2}{2} +\frac{s^2}{2}-b\,s & \mbox{ if } &  b<s\,.
	\end{array}
	\right .
\]
We will prove that the radial function
\[
    u(x) = h(|x|)= G^{-1} \left(-\log \left( \frac{|x|}{R}\right)\right)
\]
is a solution to problem \eqref{ejemplog} pointing out that, since $G^{-1}$ is discontinuous, the solution $u$ is discontinuous too.
\\
We get the radial solution
\[
    h'(r)=\frac{-1}{g\left(G^{-1}\left(-\log \left(\frac{r}{R}\right)\right)r\right)}\,,
\]
and since we take
\[
    \z(x)= \frac{-x}{|x|}\,,
\]
it is easy to prove
\[
    (\z,Du)= |Du|\,\text{ in }\, \dis (\Omega)\,,
\]
\[
    -\Div \z +g(u)^*|Du| =\frac{N}{|x|}\,\text{ as measures in }\, \Omega\,,
\]
and also
\[
    h(R) = G^{-1} (0)=0\,.
\]



\begin{thebibliography}{999}

	\bibitem{ABPP} B. Abdellaoui, L. Boccardo, I. Peral and A. Primo,
	{\it Quasilinear elliptic equations with natural growth,}
	Differential Integral Equations {\bf  20} No. 9 (2007), 1005--1020.

	\bibitem{ADS} B. Abdellaoui, A. Dall'Aglio and  S. Segura de Le\'on,
	{\it Multiplicity of solutions to elliptic  problems involving the $1$--Laplacian with a critical gradient term,}
	preprint.

	\bibitem{APP} B. Abdellaoui, I. Peral and A. Primo,
	{\it Elliptic problems with a Hardy potential and critical growth in the gradient: Non--resonance and blow--up results,}
	J. Differential Equations {\bf  239} No. 2 (2007), 386--416.

	\bibitem{Al} A. Alvino,
	{\it A limit case of the Sobolev inequality in Lorentz spaces,}
	Rend. Accad. Sci. Fis. Mat. Napoli (4) {\bf 44} (1977), 105--112.

	\bibitem{AmCM} L. Ambrosio, G. Crippa and S. Maniglia,
	{\it Traces and fine properties of a BD class of vector fields and applications,}
	Ann. Fac. Sci. Toulouse Math. (6) {\bf 14} No. 4 (2005), 527–-561.

	\bibitem{AFP} L. Ambrosio, N. Fusco and D. Pallara,
	{\it Functions of Bounded Variation and Free Discontinuity Problems,}
	Oxford Mathematical Monographs, 2000.

	\bibitem{ABCM} F. Andreu, C. Ballester, V. Caselles and J.M. Maz\'on,
	{\it The Dirichlet problem for the total variation flow,}
	J. Funct. Anal. {\bf 180} No. 2 (2001), 347--403.

	\bibitem{ACM} F. Andreu--Vaillo, V. Caselles and J.M. Maz\'on,
	{\it Parabolic quasilinear equations minimizing linear growth functionals,}
	Progress in Mathematics, 223, Birkh\"auser Verlag, Basel, 2004.

	\bibitem{An} G. Anzellotti,
	{\it Pairings between measures and bounded functions and compensated compactness,}
	Ann. Mat. Pura Appl. (4) {\bf 135} (1983), 293--318.

	\bibitem{BCN} G. Bellettini, V. Caselles and M. Novaga,
	{\it The Total Variation Flow in $\R^N$,}
	J. Differential Equations {\bf 184} (2002), 475--525.

	\bibitem{BS} C. Bennett and R. Sharpley,
	{\it Interpolation of operators,}
	Pure and Applied Mathematics, 129,  Academic Press, Inc., 1988.

	\bibitem{BMP1} L. Boccardo, F. Murat and J.P. Puel,
	{\it Existence de solutions non born\'ees pour certaines \'equations quasi--lin\'eaires,}
	Port. Math. {\bf 41} (1982), 507--534.

	\bibitem{BMP2} L. Boccardo, F. Murat and J.P. Puel,
	{\it Existence des solutions faibles des \'equations elliptiques quasi--lineaires \`a croissance quadratique,}
	Nonlinear P.D.E. and their applications, Coll\`ege de France Seminar, Vol. IV, H. Br\'ezis \& J.L. Lions (eds.), Researh Notes in Mathematics 84, Pitman, London (1983) 19--73.

	\bibitem{BMP3} L. Boccardo, F. Murat and J.P. Puel,
	{\it R\'esultats d'existence pour certains probl\`emes elliptiques quasilin\'eaires,}
	Ann. Scuola Norm. Sup. Pisa Cl. Sci. (4) {\bf 11} No. 2 (1984), 213--235.

	\bibitem{C} V. Caselles,
	{\it On the entropy conditions for some flux limited diffusion equations,}
	J. Differential Equations {\bf 250} (2011), 3311--3348.

	\bibitem{CT} M. Cicalese and C. Trombetti,
	{\it Asymptotic behaviour of solutions to $p$--Laplacian equation,}
	Asymptot. Anal. {\bf 35} (2003), 27--40.

	\bibitem{CF} G.--Q. Chen and H. Frid,
	{\it Divergence--measure fields and hyperbolic conservation laws,}
	Arch. Ration. Mech. Anal. {\bf 147} No. 2 (1999), 89--118.

	\bibitem{CTZ} G.--Q. Chen, M. Torres and W.P. Ziemer,
	{\it Gauss--Green theorem for weakly differentiable vector fields, sets of finite perimeter, and balance laws,}
	Comm. Pure Appl. Math. {\bf 62} No. 2 (2009), 242--304.

	\bibitem{D} F. Demengel,
	{\it Functions locally almost 1--harmonic,}
	Appl. Anal. {\bf 83} (2004), 865--896.

	\bibitem{EG} L.C. Evans and R.F. Gariepy,
	{\it Measure Theory and Fine Properties of Functions,}
	Studies in Advanced Mathematics, CRC Press, 1992.

	\bibitem{GMP} L. Giacomelli, S. Moll and F. Petitta,
	{\it Nonlinear Diffusion in Transparent Media: The Resolvent Equation,}
	preprint

	\bibitem{Huisken} G. Huisken and T. Ilmanen,
	{\it The Inverse Mean Curvature Flow and the Riemannian Penrose Inequality,}
	J. Differential Geom. {\bf 59} (2001), 353--438.

	\bibitem{Huisken2} G. Huisken and T. Ilmanen,
	{\it Higher regularity of the inverse mean curvature flow,}
	J. Differential Geom. {\bf 80} (2008), 433--451.

	\bibitem{H} R. Hunt,
	{\it On $L(p,q)$ spaces,}
	Enseignement Math. (2) {\bf 12} (1966), 249--276.

	\bibitem{K} B. Kawohl,
	{\it On a family of torsional creep problems,}
	J. Reine Angew. Math. {\bf 410} (1990), 1--22.

	\bibitem{KF} B. Kawohl and V. Fridman,
	{\it Isoperimetric estimates for the first eigenvalue of the $p$--Laplace operator and the Cheeger constant,}
	Comment. Math. Univ. Carolin. {\bf 44} (2003), 659--667.

	\bibitem{LL} J. Leray and J.L. Lions,
	{\it Quelques r\'esultats de Vi\v{s}ik sur les probl\`emes elliptiques non lin\'eaires par les m\'ethodes de Minty--Browder,}
	Bull. Soc. Math. France {\bf 93} (1965), 97--107.

	\bibitem{MS} J.M. Maz\'on and S. Segura de Le\'on,
	{\it The Dirichlet problem for a singular elliptic equation arising in the level set
formulation of the inverse mean curvature flow,}
	Adv. Calc. Var. {\bf 6} (2013), 123--164.

	\bibitem{MST1}  A. Mercaldo, S. Segura de Le\'on and C. Trombetti,
	{\it On the behaviour of the solutions to $p$--Laplacian equations as $p$ goes to $1$,}
	 Publ. Mat. {\bf 52} (2008), 377--411.

	\bibitem{MST2} A. Mercaldo, S. Segura de Le\'on and C. Trombetti,
	{\it On the solutions to 1--Laplacian equation with $L^1$ data,}
	J. Funct. Anal. {\bf 256} No. 8 (2009), 2387--2416.

	\bibitem{Moser} R. Moser,
	{\it The inverse mean curvature flow and $p$--harmonic functions,}
	J. Eur. Math. Soc. (JEMS) {\bf 9} (2007), 77--83.

	\bibitem{Moser2}  R. Moser,
	{\it The inverse mean curvature flow as an obstacle problem,}
	Indiana Univ. Math. J. {\bf 57} (2008), 2235--2256.

	\bibitem{S} M. \v{S}ilhav\'y,
	{\it Divergence measure fields and Cauchy's stress theorem,}
	Rend. Semin. Mat. Univ. Padova {\bf 113} (2005), 15--45.

	\bibitem{Zi} W.P. Ziemer,
	{\it Weakly differentiable functions. Sobolev spaces and functions of bounded variation,}
	Graduate Texts in Mathematics, 120, Springer--Verlag, New York, 1989.

	\bibitem{Zi2} W.P. Ziemer,
	{\it The Gauss--Green theorem for weakly differentiable vector fields,}
	Singularities in PDE and the calculus of variations, p. 233--267, CRM Proc. Lecture Notes, 44, Amer. Math. Soc., Providence, RI, 2008.

\end{thebibliography}
\end{document}